\documentclass{amse}
\usepackage{latexsym}
\makeatletter
\def\ams@pubinfo{}
\makeatother

\begin{document}
 \PageNum{1}

 \Volume{X}{X}{X}{X}
 \OnlineTime{\today}
 \DOI{00000000000000}
% \EditorNote{Received November 25, 2009, Accepted May 25, 2010}
 \EditorNote{To appear in Acta Math. Sin. (Engl. Ser.), DOI:10.1007/s10114-010-9696-9.}

 \abovedisplayskip 6pt plus 2pt minus 2pt \belowdisplayskip 6pt plus 2pt minus 2pt
%%%%%%%%%%%%%%%%%%%%%%%%%%%%%%%%%%%%%%%%%%%%%%%%%%%%%%%%%
%%-------------------       the beginning of  Author's Definitions       -------------------%%
%%   Note:  The following are some samples, you may replace it with your own definitions.
%%          Please leave it EMPTY if you don't have any definitions.
%%
\numberwithin{equation}{section}
\numberwithin{example}{section}
\numberwithin{rem}{section}
%%-------------------         the end of  Author's Definitions           -------------------%%

 \AuthorMark{Thomas Kaijser}                             %%%  appear on the head of even pages  %%%

 \TitleMark{On Markov chains induced by partitioned %transition probability
   matrices}                  %%%  appear on the head of odd pages  %%%

\title{On Markov chains induced by partitioned transition probability
matrices                        %%%   Main Title of your paper  %%%
\footnote{  %Supported by National Natural Science Foundation of China (Grant No. )
}}                                                  %%%   the Fund which you are supported by  %%%

\author{Thomas  \uppercase{Kaijser}}                                              %%%  1st Author's information   %%%
{ Department of Mathematics
 and  Information Coding Group,\\
 Link\"{o}ping University,
 S-581 83 Link\"{o}ping, Sweden.\\ 
 E-mail\,$:$ thkai@mai.liu.se }
\maketitle%

\Abstract{Let \(S\) be a denumerable state space and let \(P\) be a
  transition probability matrix on \(S\). If a denumerable set \({\cal
    M}\) of nonnegative matrices is such that the sum of the matrices
  is equal to \(P\), then we call \({\cal M}\) a {\em partition} of
  \(P\).\newline
  Let \(K\) denote the set of probability vectors on \(S\).  To every
  partition \({\cal M}\) of \(P\) we can associate a transition
  probability function \({\bf P}_{\cal M}\) on \(K\) defined in such a
  way that if \(p \in K\) and \(M \in {\cal M}\) are such that
  \(||pM|| > 0 \), then, with probability \(||pM||\), the vector \(p\)
  is transferred to the vector \(pM/||pM||\).  Here \(||\cdot||\)
  denotes the \(l_1-norm.\)\newline
  In this paper we investigate convergence in distribution for Markov
  chains generated by transition probability functions induced by
  partitions of transition probability matrices.  The main motivation
  for this investigation is the application of the convergence results
  obtained to filtering processes of partially observed Markov chains
  with denumerable state space.  }  % the abstract

        % the keywords

\MRSubClass{Primary 60J05; Secondary
  60F05.}  % MR(2000) Subject Classification

\section{Introduction}
\subsection{ The filtering process.} 
Let \(S\) denote a denumerable
set, \(\{X_n, n=0,1,2,...\}\) be an  aperiodic, 
positively recurrent  Markov chain with  
\(S\) as state space,
let \(A\) denote an ``observation space'',
let  \(g:S \rightarrow A\) denote a "lumping" function of 
the state space \(S\) and define 
\(Y_n = g(X_n)\). Let \(Z_n\) denote 
the conditional  distribution of \(X_n\) given \(Y_1,Y_2...,Y_n.\)
  
In this paper we shall present sufficient conditions both 1) for when 
the distributions of the {\em filtering process} \(\{Z_n, n=1,2,...\}\)
{\em do } converge in distribution to a 
unique limit distribution independent 
of the initial distribution, and 2) when they {\em do not}.

Our paper is centered around a notion we call a 
{\em  partition of a transition
  probability
matrix} (tr.pr.m) and we shall soon see how this notion connects with the
filtering process just described.

\subsection{Partitions of  transition probability matrices.}
Let \(S\) be a denumerable set. The set of all tr.pr.ms on \(S\)
will be denoted \(PM(S)\). 
A denumerable set 
\[{\cal M}=\{M(w): w \in {\cal
  W}\}\]
of 
nonnegative
\(S\times S\) matrices such that 
\[\sum_{w \in {\cal W}}M(w) = P \] 
will 
be called
a {\bf partition of \(P\)}. 

The set of all partitions of \(P\) will be denoted \(G(S,P)\).
We define 
\[
G(S) = \cup_{P \in PM(S)}\; G(S,P)
\]
and call an element in \(G(S)\) simply a {\bf partition}.

We denote the \(i,j-th\) element of a matrix \(M\) by 
\((M)_{i,j}.\)
\begin{example}\label{ex1.1} Let \(S\) be a denumerable set, let
\(P \in PM(S)\), let \(A\) be another denumerable set
and let \(g:S \rightarrow A\) be a ``lumping'' function
from \(S\) to \(A\). For each element \(a \in A\)
 we define the matrix 
M(a) by
\[
(M(a))_{i,j} = (P)_{i,j}\;\;  if \;\; g(j)= a 
\]
\[
(M(a))_{i,j} = 0\;\;  if \;\; g(j)\not = a. 
\]
\end{example}
Obviously \(\{M(a):a \in A\}\) is a partition of \(P\).
A partition defined in this way by a ``lumping'' function
\(g\) is simply called the {\em partition determined by the ``lumping''
function \(g\)}.

\begin{example} \label{ex1.2}
Let \(S\) be a denumerable set, let \(A\) be another denumerable set,
let \(P \in PM(S)\) and let \(R\) be a tr.pr.m from \(S\) to \(A\).
For each element \(a \in A\) we define the matrix 
\(M(a)\) by
\[
(M(a))_{i,j} = (P)_{i,j}(R)_{j,a}.\;\;  
\]
\end{example}
Again it is easily seen that \(\{M(a):a \in A\}\) is a partition of 
\(P\).
A partition defined in this way by an {\em observation matrix}
\(R\) is simply called the {\em partition determined by the
  observation matrix \(R\)}.
\newline
\begin{rem} Note that Example \ref{ex1.1} is
a special case of Example \ref{ex1.2} since if \(g\) is a ``lumping'' function 
we can define a tr.pr.m \(R\) from \(S\) to \(A\) simply by
\[(R)_{j,a} =1 \;\;if \;\; g(j)=a\]
\[(R)_{j,a} =0 \;\;otherwise.\]
\end{rem}
\subsection{ The transition probability function 
\({\bf P}_{{\cal M}}\).}  
Let again \(S\) denote a denumerable set. The set of all 
probability vectors on \(S\)
will be denoted by \(K\).  
Thus
\begin{equation}\label{eq1.1.1}
K=\{x= ((x)_i, i \in S) : (x)_i \geq 0, \; \sum (x)_i = 1\}.
\end{equation}

We consider elements in \(K\) as row vectors. We denote the \(i-th\)
coordinate
of a vector \(x\) by \((x)_i\).
We let \(||\cdot||\) denote the \(l_1-norm\) and 
introduce a distance function \(\delta\) on \(K\) by using the 
\(l_1-norm\). Thus 
\[
\delta(x,y) = ||x-y|| = \sum_{i \in S} |(x)_i - (y)_i|.
\]
We let \({\cal E}\) denote the Borel field generated by the
metric \(\delta\) and let \({\cal P}(K)\) denote the set of all
probability measures on \((K, {\cal E}).\)

Now, to every partition \({\cal M} =\{M(w): w \in {\cal W}\}\in G(S)\),
we  can define a {\bf transition probability function}
(tr.pr.f) \({\bf P}_{{\cal M}}: K\times {\cal E} \rightarrow [0,1]\) 
on \((K,{\cal E})\)
by
\begin{equation}\label{Ny2eq1.1}   
{\bf P}_{{\cal M}}(x,B) = 
\sum_{w \in {\cal W}_{{\cal M}}(x,B)}
||xM(w)||, \;\;
\end{equation}
where 
\begin{equation}\label{eq1.5.5}
{\cal W}_{{\cal M}}(x,B) =\{w\in {\cal W}: ||xM(w)||>0, 
xM(w)/||xM(w)|| \in B\}.
\end{equation}
That \({\bf P}_{{\cal M}}\), as defined by (\ref{Ny2eq1.1}) 
and (\ref{eq1.5.5}), does indeed define a tr.pr.f is not 
very difficult to 
prove. 
We sketch a proof at the end of the next section.

Next, let \({\bf P}_{{\cal M}}^n(\cdot,\cdot)\) denote the \(n\)-step
tr.pr.f  defined recursively by
\[
   {\bf P}_{{\cal M}}^1(x,B) =
   {\bf P}_{{\cal M}}(x,B), \;\; \;\; x \in K, \;\; B \in {\cal E},
\]
\[
   {\bf P}_{{\cal M}}^{n+1}(x,B) =
\int_K    {\bf P}_{{\cal M}}^n(y,B) {\bf P}_{{\cal M}}(x,dy),
\;\;\; x \in K, \;\; B \in {\cal E}, \;\;n=1,2,...\;.
\]
Let \(C[K]\) denote the set of all real, continuous, bounded, 
functions on \(K\). If the tr.pr.f 
\({\bf P}_{{\cal M}}(\cdot,\cdot)\) is such that there exists
a probability measure \(\mu \in {\cal P}(K)\)
such that  
\[
\lim_{n \rightarrow \infty} 
\int_K u(y){\bf P}_{{\cal M}}^n(x,dy) =  \int_K u(y)\mu(dy),
\;\forall u \in C[K], \;\forall x \in K\]
then we say that \(
{\bf P}_{{\cal M}}(\cdot, \cdot)
\) is {\bf asymptotically stable}.

The main purpose of this paper is to give a {\em sufficient condition}
 for asymptotic
stability of \({\bf P}_{{\cal M}}(\cdot, \cdot)\) when 
the tr.pr.m \(P\) on \(S\) is irreducible,
aperiodic and positively recurrent. Another purpose is to give 
a {\em sufficient condition} for when 
\({\bf P}_{{\cal M}}(\cdot, \cdot)\) is {\em not} asymptotically stable. 

\subsection{ The interrelationship with the filtering
  process.}
Let again \(S\) and \(A\) be denumerable sets. 
 Let \(\{(X_n,Y_n), n=0,1,2,...\}\) be  a 
{\em hidden Markov chain} determined by a tr.pr.m \(P \in PM(S)\),
a tr.pr.m \(R\) from \(S\) to \(A\), and an initial distribution \(p\).
Define
\[
Z_{n,i} =Pr[X_n=i| Y_1,Y_2,...,Y_n], \; i \in S,
\;\;n =1,2,...
\]
set 
\[
Z_n=(Z_{n,i}, i \in S)
\]
and set
\[
\mu_{n,p} = probability \;\; distribution \;\; of\;\; Z_n.
\]
Let \({\cal M}\) be the partition of Example 1.2 and let 
\({\bf P}_{{\cal M}}\) be the tr.pr.f induced by 
\({\cal M}\).
Then, for \(n=1,2,...\) and \(p \in K\)
\begin{equation}\label{C3.2}
\mu_{n,p}(\cdot) =  {\bf P}_{{\cal M}}^n(p,\cdot).
\end{equation}
\begin{rem} That (\ref{C3.2}) is true is easy to prove. (For the
case when the observation matrix is determined by a ``lumping''
function
see e.g \cite{Bla57},\cite{Ast65}, \cite{Rud73}.)
\end{rem}
\begin{rem} 
The stochastic quantity \(Z_n\) as defined above is
often called
the {\em conditional state distribution} (at time \(n\)), and the sequence
\(\{Z_n, n=1,2,...\}\) 
is  often called the {\em filtering process}.
\end{rem}

For some general basic theory regarding hidden Markov chains
see e.g \cite{RJ86} and \cite{EM02}.

\subsection{Previous work.}
In the  classical paper \cite{Bla57} by D. Blackwell
from 1957 the author proves that 
the tr.pr.f \({\bf P}_{{\cal M}}\) has a unique invariant measure
if \(S\) is finite,
the partition is determined by a ``lumping'' function and the tr.pr.m
\(P \in PM(S)\) has  ``nearly identical rows and no element 
which is very small''.

In the paper \cite{Kai75} from 1975 asymptotic stability was proven  
when \(S\) is finite, the partition 
\({\cal M}= \{M(w): w \in {\cal W}\}\) is determined by a
``lumping'' function, the associated tr.pr.m \(P\) is aperiodic and 
irreducible,  and 
the following condition is satisfied:
\newline

{\bf Condition A.} {\em There exists an integer \(N\) and a sequence
\(w_1, w_2,...,w_N\) of elements in \({\cal W}\)
such that
the matrix product 
\(M(w_1)M(w_2)...M(w_m)\) is a nonzero matrix 
with the property that
if  
\((M(w_1)M(w_2)...M(w_m))_{i_1,j_1} \not = 0 \) 
and also
\((M(w_1)M(w_2)...M(w_m))_{i_2,j_2} \not = 0 \) 
then also
\((M(w_1)M(w_2)...M(w_m))_{i_1,j_2} \not = 0 \) 
and 
\((M(w_1)M(w_2)...M(w_m))_{i_2,j_1} \not = 0. \) 
}
\newline

Next, let again \(S\) be finite, 
and \({\cal M}=\{ M(w) : w \in {\cal W}\}\in G(S)\). 
For \(n=1,2,...\) define
\[
{\cal M}^n = \{M(w_1)M(w_2)...M(w_n) : w_i \in {\cal W}, \;
1 \leq i \leq n\},
\]
set
\[
{\cal M}^* =\cup_{n=1}^\infty {\cal M}^n, 
\]
set
\[{\cal C} = \{\alpha M : \alpha >0, M\in {\cal M}^*
\}\]
and 
let the set 
\(\;{\overline{\cal C}}\)
be defined as  the closure of 
\(\;{\cal C}\) under the usual topology in \({\tt R}^{S\times S}\)
where \({\tt R}\) denotes the set of real numbers.

The following condition, which we call
{\bf Condition KR} was introduced by F. Kochman and J. Reeds 
in the paper \cite{KR06} from 2006.
\newline

{\bf Condition KR.} {\em The set  \(\overline{{\cal C}}\) contains a
matrix of {\bf rank 1}.}
\newline

In \cite{KR06} it is proved that if \(S\) is finite, 
the tr.pr.m \(P\) is aperiodic
and irreducible, the partition \({\cal M}\) of \(P\) is determined
by an observation matrix and Condition KR holds then 
the induced tr.pr.f
\({\bf P}_{\cal M}\) is asymptotically stable.

\subsection{The main theorem.}
Let \(S\) be a denumerable set. 
If  \(M\) is an  \(S\times S\) matrix
we define the norm \(||M||\) by
%\begin{equation}
%\label{i1eq11}
\[
||M|| = \sup \{||xM||: ||x|| =1, x \in {\tt R}^S\}.
%\end{equation}
\]
Next let \({\cal U}\)  denote the set of \(S-dimensional\) vectors 
specified
by
\[
{\cal U}=\{u=((u)_i, i \in S) : u_i\geq 0, \;\;
\;\;and \; \; 
\sup\{ u_i: i \in S \} =\;1 \;\;\},
\]
and 
let  \({\tt W}\) denote  the set of 
\(S\times S\)  matrices  defined by
\[
{\tt W} =\{ W=u^c v : u \in {\cal U} , v \in K \}
\]
where  \(u^c\) denotes the transpose of \(u\).
Note that if \(W \in {\tt W}\),
then \(||W||=1\) since \(0\leq \sum_{j\in S}u_iv_j = u_i\leq 1\)
for all \(i \in S\) and  
\(\sup_i u_i =1. \)
We call an element in \({\tt W}\) a {\em nonnegative rank 1 matrix of
  norm 1}.

Let \({\cal M} =\{M(w):w \in {\cal W}\}\) denote a partition.
If \(w_1,w_2,...,w_m\) is a finite 
sequence of elements in  \({\cal W}\) we use the notations
\[
{\bf w^m}=(w_1,w_2,...w_m),
\]
and
\[
{\bf M}({\bf w^m}) = M(w_1)M(w_2)...M(w_m).
\]

We shall next introduce two 
conditions for a partition
\({\cal M}=\{M(w) \,: \, w \in {\cal W}\} \in G(S)\). 
Our first condition is a rather straight forward generalization
of Condition KR.
Here and throughout this paper we let 
\(e^i, \; i \in S\) denote the vector in \(K\) defined by
%\begin{equation}
\[
(e^i)_i=1.
\]
\[\]
{\bf Condition B1.}
{\em There exists 
a nonnegative rank 1 matrix 
\(W = u^cv\) of norm 1, 
a sequence of integers \(\{n_1,n_2,...\}\) and a sequence
\({\bf w_j^{n_j}}
, j=1,2,...\)
of sequences  
\({\bf w_j^{n_j}}=
w_{1,j},w_{2,j},...,w_{n_j,j}, \;w_{k,j}\in {\cal W}, 
1\leq k \leq n_j\), 
such that 
\(||{\bf M}({\bf w_j^{n_j}})||>0,\; \;\;j =1,2,...\) 
and such that for all 
\(i\in S\) 
\[
\lim_{j\rightarrow \infty} 
||e^{i}{\bf M}({\bf w_j^{n_j}})/
||{\bf M}({\bf w_j^{n_j}})|| - 
 e^{i}W|| = 0.
\]
}

It is not difficult to prove that if in Condition B1 the underlying
set \(S\) is finite
then Condition B1 is equivalent to Condition KR.

In order to define our next condition we first need to introduce the
well-known notion {\bf barycenter}. 
The barycenter of a measure 
\(\mu \in {\cal P}(K)\) is defined as that vector 
\({\overline b}(\mu)\in K\) whose \(i-th\) coordinate 
\(({\overline b}(\mu))_i\) 
is
defined by
\begin{equation}\label{i1eq12}
({\overline b}(\mu))_i = \int_K (x)_i \mu(dx).
\end{equation}
That the vector \({\overline b}(\mu)\) belongs to \(K\) follows 
immediately from
the fact that \(\sum_{i\in S}(x)_i=1\).
We let \({\cal P}(K|q)\) denote the subset  
of \({\cal P}(K)\) 
such that 
each \(\mu \in  
{\cal P}(K|q)\) has 
barycenter 
equal to \(q\).

We are now ready to introduce the condition under which the main
theorem of this paper is proved. Let 
\(S\) be a denumerable set,
let \(P \in PM(S)\)
be irreducible,
aperiodic
and positively recurrent, 
let \(\pi\) denote the unique 
probability vector in \(K\) such that \(\pi P=\pi\) and 
let \({\cal M}=\{M(w) \,: \, w \in {\cal W}\}\) 
be a partition of  \(P\). 
\newline

{\bf Condition B.} 
 {\em For every \(\rho > 0\) there exists an element 
\(i_0 \in S\) such that if \(C\subset K\) is a compact set 
satisfying
\begin{equation}\label{I1eq21}
\mu(C\cap \{x:(x)_{i_0}\geq (\pi)_{i_0}/2 \})  \geq (\pi)_{i_0}/3,
\;\;\;\;\forall \; \mu \in {\cal P}(K|\pi), 
\end{equation}
then we can find 
an integer \(N\), 
and a sequence \(w_1,w_2,...w_N\) 
of elements in \({\cal W}\),
such that, if we set
\[
{\bf M}({\bf w^N}) = M(w_1)M(w_2)...M(w_N),
\]
then 
\[
||e^{i_0}{\bf M}({\bf w^N})|| > 0
\]
and if \(x \in C\cap \{x:(x)_{i_0}\geq (\pi)_{i_0}/2 \}\) then also 
\[
||(x{\bf M}({\bf w^{N}})/||x{\bf M}({\bf w^{N}})|| -
e^{i_0}{\bf M}({\bf w^{N}})/||e^{i_0}{\bf M}({\bf w^{N}})||)|| < \rho.\]
}
\newline
That there exists a compact set \(C\) such that 
(\ref{I1eq21}) holds is proved in section 4.

It is not very difficult to prove that Condition B1 implies
Condition B when \(P\) is aperiodic, irreducible and positively
recurrent, a fact we shall prove in section 9.

The main theorem of this paper reads as follows:
\begin{theorem}\label{thm1.1} 
%{\em
Let 
\(S\) be a denumerable set, let \(P \in PM(S)\)
be irreducible, aperiodic
and positively recurrent, 
let \(\pi \in K\) satisfy \(\pi P =\pi\),  
let 
\({\cal M}\) be a partition of \(P\) and 
let \({\bf P}_{{\cal M}}\) be the tr.pr.f induced by
\({\cal M}\). Suppose also that
Condition B holds. Then 
\({\bf P}_{{\cal M}}\) is asymptotically stable.
%} 
\end{theorem}

\subsection {Exceptional cases.} 
One consequence of asymptotic stability is that there only exists
one {\em invariant measure}. Therefore, if \({\bf P}_{{\cal M}}\)
fulfills the hypotheses of Theorem 1.1 
then the equation
\begin{equation}\label{I1eq1.6.1}
\int_{K} {\bf P}_{{\cal M}}(dx, B)\mu(dx) = \mu(B), \;\;\;
\forall B \in {\cal E}
\end{equation}
has a unique solution in \({\cal P}(K)\).
In the paper \cite{Bla57} D. Blackwell conjectured 
that the equation (\ref{I1eq1.6.1}) has a unique solution 
if \(S\) is finite, \(P \in PM(S)\) is indecomposable
and the partition is determined by a "lumping" function on \(S\).
However, there are counterexamples to 
this conjecture and one such counterexample was presented in
\cite{Kai75}. 
In fact, already in 1974, H. Kesten constructed an 
example, not published before, which shows that the tr.pr.f 
\(\,{\bf P}_{{\cal M}}\) 
can  in fact even be
{\bf periodic} (\cite{Kes74}). 
 In section 11 we present this counterexample. 
 
In 
section 11
we also state and 
prove a theorem with hypotheses that
guarantee that 
\({\bf P}_{\cal M}\)
is {\bf not}
asymptotically stable and describe a whole class of tr.pr.ms and
partitions such that the induced tr.pr.f is {\bf not}
asymptotically stable.

\subsection{The plan of the paper.} 
In section 2 we introduce some further notations and concepts.
In section 3 we present a few basic lemmas for
tr.pr.fs induced by partitioned tr.pr.ms.
In section 4 we introduce the
well-known
Kantorovich distance and in sections 5 and 6 we prove some results 
for probability measures with equal barycenter. One result is that
the set \(P(K|q)\) is a {\bf tight} set for every \(q \in K\). 

In section 7 we introduce a property which we call the 
{\bf shrinking property} for tr.pr.fs on metric spaces, and prove
an auxiliary ergodic theorem for Markov chains induced by such tr.pr.fs.
In section 8 we prove the main theorem of this paper 
(Theorem \ref{thm1.1})
by verifying that
the hypotheses of the auxiliary theorem of section 7 is fulfilled.
In section 9  we verify that Condition B1 implies Condition B,
and in section 10 we present two  random walk examples satisfying
Condition B1.

In section 11, as mentioned above, we consider exceptional cases.
In section 12 we consider convex functions and 
state some inequalities reminiscent of inequalities obtained
by H. Kunita  in his classical paper 
\cite{Kun71} from 1971,
and in section 13 finally, 
we generalize Blackwell's entropy formula for the entropy rate
of functions of finite-state Markov chains presented in \cite{Bla57}
to Markov chains with denumerable state space.
By using  convexity properties proved in section 12
we can also give lower and upper bounds for the entropy rate.

\section{ Some further notations and concepts.}
Let \(S\) be a denumerable set and \(K\) the set defined by
 (\ref{eq1.1.1}). 
We let \(PM_{e}(S)\) denote the set of tr.pr.ms on \(S\)
such that  if \(P \in PM_{e}(S)\) then there exists a unique 
vector 
\(\pi \in K\) with {\em positive} coordinates such that 
\begin{equation}\label{eq1.2.1}
\pi P = \pi,
\end{equation}
and we let 
\(PM_{ae}(S)\) denote the subset of
\(PM_{e}(S)\) 
which consists of tr.pr.ms 
 which are {\em aperiodic}.
Recall from the general theory on Markov chains that if
\(P \in PM(S)\) is aperiodic, irreducible, and positively recurrent
then \(P \in PM_{ae}(S)\) and if also \(\pi\) satisfies
(\ref{eq1.2.1}) then
\begin{equation}\label{eq2.1.1}
\lim_{n \rightarrow \infty} ||xP^n - \pi|| = 0, \;\; \forall x \in K.
\end{equation}
(See e.g \cite{Lin92}, Chapter 2.) 

For \(u \in C[K]\) we define \(\gamma(u)\) by
\[
\gamma(u) = \sup \{ |u(x)-u(y)|/||x-y|| : x,y \in K, \;x\not = y \},
\]
we define \(Lip[K] = \{u \in C[K] : \gamma(u) < \infty\}\)
and \(Lip_1[K]=\{u \in Lip[K] : \gamma(u) \leq  1\}\).
We let \(F[K]\) denote the set of real, bounded functions on \(K\), and we
let \(B[K]\) denote the set of real \({\cal E}-measurable\) 
functions on \(K\). For \(u \in F[K]\)
we define \(||u|| = \sup_{x \in K} |u(x)|,\)
we define
\[osc(u) = \sup \{u(x)-u(y): x,y \in K\}\]
and if \(A \subset K\) we define
\begin{equation}\label{eq2.15}
osc_A(u) =\sup \{u(x)-u(y): x,y \in A\}.
\end{equation}
We set \(l_1(S) = \{x \in {\tt R}^S : 
\sum_{i\in S} |(x)_i|< \infty \}\) and 
\(l^{\infty}(S) = \{x \in {\tt R}^S : 
\sup \{|(x)_i|:i \in S\} < \infty\}   \).

Next, we present some trivial facts regarding partitions, facts
which we state without proof. 
\newline
1) Let
\({\cal M}_1=\{A(w_1): w_1 \in {\cal W}_1\} \in G(S,P_1)\),
\(\;{\cal M}_2=\{B(w_2): w_2 \in {\cal W}_2\} \in G(S,P_2)\) and 
define
\(
{\cal M}_3=\{A(w_1)B(w_2): 
w_1 \in {\cal W}_1,\;w_2 \in {\cal W}_2\}.
\)
Then 
\[
{\cal M}_3  
 \in G(S,P_1P_2).
\]
We write \(\;{\cal M}_3 ={\cal M}_1 {\cal M}_2 \;\)
and call 
\({\cal M}_3\) the {\em product} of
\({\cal M}_1\) and \( {\cal M}_2\).
\newline
2) Let \({\cal M}_i\in G(S), \;i=1,2,3\). Then 
\[
({\cal M}_1{\cal M}_2){\cal M}_3 =
{\cal M}_1 ({\cal M}_2{\cal M}_3).
\]
Therefore if 
\({\cal M}_i \in G(S), i=1,2,...,n\)
then 
\(
{\cal M}_1{\cal M}_2...{\cal M}_n \)
is well-defined.
If \({\cal M}_i= \{M_i(w_i): w_i \in {\cal W}_i\} \in G(S),
\;i=1,2,...n\) 
we write
\[
{\bf w}^n  = (w_1,w_2,...w_n)
\;\;and\;\;
{\cal W}^n =\{(w_1,w_2,...,w_n):\; w_i \in {\cal W}_i\}\]
and write
\[
 {\cal M}^n =
{\cal M}_1{\cal M}_2...{\cal M}_n =
\{M_1(w_1)M_2(w_2)...M_n(w_n): \;w_i \in {\cal W}_i\}.
\]
We often denote an element in
\( {\cal M}^n\) by \({\bf M}({\bf w}^n)\).

 From above follows also that 
\begin{equation}\label{eq1.4.12}
if \;\;{\cal M}_i= {\cal M}\in G(S,P), \; i=1,2,...n,\;\; 
then \;\;
 {\cal M}^n \in G(S,P^n).
\end{equation}

Next, let \({\cal M}=
\{M(w): w \in {\cal W}\} \in G(S)\) and 
let \({\bf P}_{{\cal M}}\) be the tr.pr.f
on \((K,{\cal E})\) induced by \({\cal M}\).
The tr.pr.f \({\bf P}_{{\cal M}}\) determines two 
mappings, one which is defined on the set \(F[K]\) of real bounded
functions, and one which is defined on the set \({\cal P}(K)\) 
of probability measures. Thus we define 
\(T_{{\cal M}}: F[K] \rightarrow F[K] \)
by
\[
T_{{\cal M}}u(x)=
\sum_{w \in {\cal W}_{{\cal M}}(x)}u(xM(w)/||xM(w)||)||xM(w)||
\]
where \({\cal W}_{{\cal M}}(x)\) is defined by 
\begin{equation}\label{eq1.5.1}
{\cal W}_{{\cal M}}(x)
=\{w \in {\cal W} : ||xM(w)||> 0\},
\end{equation}
and we define \({\breve P}_{{\cal M}}: {\cal P}(K) \rightarrow 
{\cal P}(K)\)
by
\[
{\breve P}_{{\cal M}}\mu(B)=
\int_{K} {\bf P}_{{\cal M}}(x,B)\mu(dx), \; B \in {\cal E}.
\]
We call \(T_{{\cal M}}\) the {\bf transition operator} induced by
\({\cal M}\) and we call
\( {\breve P}_{{\cal M}}\)
the  {\bf transition probability operator} 
induced by \({\cal M}\).

For \(u \in B[K]\) and \(\mu \in {\cal P}(K)\) we shall - when
convenient -
write
\[
\langle u,\mu \rangle = \int_K u(x)\mu(dx).
\] 
It is well-known that 
\begin{equation}\label{I1eq1.5.3}
\langle T_{{\cal M}}u,\mu\rangle =\langle u, 
{\breve P}_{{\cal M}}\mu\rangle, \;\;u\in B[K], \;\mu \in {\cal P}(K).
\end{equation}
(See e.g. \cite{Rev75}, chapter 1, section 1.)

We end this section by sketching a proof of the
fact that the function
\({\bf P}_{{\cal M}}: K\times {\cal E} \rightarrow [0,1]\) 
defined by
(\ref{Ny2eq1.1}) and   
(\ref{eq1.5.5}) does indeed determine a tr.pr.f.

Firstly,
that \({\bf P}_{{\cal M}}(x, \cdot)\) for each \(x\in K\)
satisfies the relations that
define a probability measure is easy to verify from the definition
of \({\bf P}_{{\cal M}}\).

Secondly in order to prove that
 \({\bf P}_{{\cal M}}(\cdot, B)\) is \({\cal E}-measurable\),
we first note that it suffices to prove
that for each  \(w \in {\cal W}\) 
and each \(B \in {\cal E}\)  
the function  
\(u_{B}:K \rightarrow {\tt R}\) defined by 
\[ u_{B}(x) =||xM(w)||, \;if\; ||xM(w)||>0 \;
and \; xM(w)/||xM(w)|| \in B,\]
 \[ u_{B}(x) = 0, \; otherwise
\]
is \({\cal E}-measurable\). To do this let  
\(
{\cal F} = \{B \subset K : u_B \;\;is \;\;{\cal E}-measurable\}.
\)
From the definition of \(u_B\) it is easy to verify 
that 
\({\cal F}\) is in fact a \(\sigma-algebra\). 
 Therefore it remains to
show that the function \(u_B(x)\) is \({\cal E}-measurable\)
if \(B\) is an {\em open} set. But this is easily
done by using very standard type of arguments.

\section{ Some basic lemmas.}
In this section we shall collect a few basic  equalities and
inequalities which will be used in the proof of Theorem 1.1.

We start with the following trivial {\em scaling property} 
for matrix products
which is used frequently. We omit the proof. 
\begin{lemma}\label{lem3.1} Let \(A\) and \(B\) denote two matrices
(not necessarily of finite dimension) and assume that
\(
AB\) is well defined. Let \(x\) be a row vector
and assume also that 
\newline
1) \(xA\) is well defined, 2)  
\(\;\; 0 < ||xA|| < \infty\;\) and 3)  
\(\;\; 0 < ||xAB|| < \infty\). Then
\begin{equation}\label{G2eq1} 
xAB/||xAB||= (xA/||xA||)B/||(xA/||xA||)B||.
\end{equation}
\end{lemma}

Using the scaling property and 
the fact that 
\begin{equation}\label{eq1.3.2}
\sum_{w \in {\cal W}} ||xM(w)|| = ||xP|| = 1, \;\;\forall x \in K
\end{equation}
 it follows 
easily that
if 
\({\cal M}_1=\{M_1(w_1):w_1 \in {\cal W}_1\} \in {\cal G}(S)\)
and
\({\cal M}_2=\{M_2(w_2):w_2 \in {\cal W}_2\} \in {\cal G}(S)\)
then
\begin{equation}\label{G2eq11}
\sum_{w_2 \in {\cal W}_2}||xM_1(w_1)M_2(w_2)|| = ||xM_1(w_1)||.
\end{equation}
a relation which is useful when proving the following
lemma.
\begin{lemma}\label{lem3.2}
Let
\({\cal M}_1=\{M_1(w_1):w_1 \in {\cal W}_1\} \in {\cal G}(S)\)
and
\({\cal M}_2=\{M_2(w_2):w_2 \in {\cal W}_2\} \in {\cal G}(S).\)
Then
\[
T_{{\cal M}_1{\cal M}_2} =
T_{{\cal M}_1}T_{{\cal M}_2}. 
\]
\end{lemma}
{\em Proof}.  
Let \(u \in F[K]\) and \(x \in K\). 
From (\ref{G2eq1}) and  
(\ref{G2eq11}) follows
\(
T_{{\cal M}_1{\cal M}_2}u(x) =\)
\[ 
\sum_{(w_1,w_2) 
\in {\cal W}_{{\cal M}^2}^2(x)}
u(xM_1(w_1)M_2(w_2)/||xM_1(w_1)M_2(w_2)||)\cdot
||xM_1(w_1)M_2(w_2)|| = \]
\[
\sum_{w_1\in {\cal W}_{{\cal M}_1}(x)} 
T_{{\cal M}_2}u(xM_1(w_1)/||xM_1(w_1)||)(||xM_1(w_1)||) = 
T_{{\cal M}_1}T_{{\cal M}_2}u(x)
\]
from which the lemma follows. \(\;\Box \)
\begin{corollary}\label{cor3.3}
Let 
\({\cal M}_1, \;{\cal M}_2 \;\in {\cal G}(S)\).
Then

\[
{\breve P}_{{\cal M}_2} {\breve P}_{{\cal M}_1} =
{\breve P}_{{\cal M}_1 {\cal M}_2}.
\]
\end{corollary}

\begin{corollary}\label{cor3.4}
 Let 
\({\cal M} \in {\cal G}(S)\).
Then for \(n=2,3,...\)
\[
T^n_{{\cal M}}
 =
T_{{\cal M}^n}  \;\;and \;\;
{\breve P}^n_{{\cal M}}
 =
{\breve P}_{{\cal M}^n}  \;\; 
\]
\end{corollary}

The following {\em universal inequality} 
was in principal already proved in 
\cite{Kai75},  section 4. 
 
\begin{lemma}\label{lem3.5}
Let \(S\) be a denumerable set 
and suppose \({\cal M}
\in {\cal G}(S)\).
Then, 
\begin{equation}\label{G3eq0}
u \in Lip[K] \;\Rightarrow T_{{\cal M}}u \in Lip[K] 
\end{equation}
and   
\begin{equation}\label{R3eq1} 
\gamma(T_{{\cal M}}u) \leq 3 \gamma(u). 
\end{equation}
\end {lemma}
{\em Proof.} Since the arguments needed can be found in \cite{Kai75}
we will be rather brief. 
We shall first prove that 
\begin{equation}\label{C5.eq0}
\gamma(T_{{\cal M}}u) \leq ||u||+2 \gamma(u) , \;\;\forall  u \in Lip[K]. 
\end{equation}
Let \(x,y \in K\) and \(u \in Lip[K]\).
We want to estimate
\(|T_{{\cal M}}u(x)-T_{{\cal M}}u(y)| \).
Let
\( S_1= \{i: (x)_i>0, \;(y)_i> 0\}, \)
\( S_2= \{i: (x)_i>0, \;(y)_i= 0\}, \)
\( S_3= 
\{i: (x)_i=0, \;(y)_i> 0\}, \)
\({\cal M}=\{M(w) : w \in {\cal W}\}\) 
and  for \(i \in S\) set 
\[ {\cal W}_i=\{w \in {\cal W} : ||e^iM(w) || >0\}.\]

Using the fact that \(T_{{\cal M}}u(x)\) can be represented by
  \[
T_{{\cal M}}u(x)=
\;\sum_{i \in S_1} (x)_i \sum_{w \in {\cal W}_i}
u(xM(w)/||xM(w)||)\cdot||e^iM(w)|| +\]
\[
\sum_{i \in S_2} (x)_i \sum_{w \in {\cal W}_i}
u(xM(w)/||xM(w)||)\cdot||e^iM(w)|| \]
and using a similar representation for 
\(T_{{\cal M}}u(y)\) 
(with \(S_2\) replaced by \(S_3\))
it is easily proved, by  using the triangle
inequality,
that
\[|T_{{\cal M}}u(x)-T_{{\cal M}}u(y)| \leq \]
\[
||u||\cdot ||x-y||
+\]
\begin{equation}\label{G3eq11}
\gamma(u) \sum_{i \in S_1} (y)_i \sum_{w \in {\cal W}_i}
  ||(xM(w)/||xM(w)||-yM(w)/||yM(w)||)||\cdot||e^iM(w)||.
\end{equation}

In order to prove that the last term is less
than \(2\gamma(u)||x-y||\) 
we shall use the following 
inequality.  
{\em Let \(a\) and \(b\) be two nonzero
vectors 
in a normed
vector space. Then
\begin{equation}\label{R3eq6}
||(a/||a|| -b/||b||)|| \leq 2||a-b||/||b||.
\end{equation}
}
Using (\ref{R3eq6})
 we find
that 
\begin{equation}\label{F3eq1}
||(xM(w)/||xM(w)||-yM(w)/||yM(w)||)||\leq 2 ||xM(w)-yM(w)||/||yM(w)||
\end{equation}
if \(||xM(w)||\cdot ||yM(w)|| > 0, \) 
and by using (\ref{F3eq1}), the triangle inequality and change of
summation order, we obtain  that 
\[
\sum_{i \in S_1} (y)_i\sum_{w \in {\cal W}_i}
||(xM(w)/||xM(w)||-yM(w)/||yM(w)||)||\cdot||e^iM(w)||
\leq
\]
\[
2\sum_{w \in {\cal W}}||xM(w)-yM(w)|| 
\sum_{i \in S} (y)_i||e^iM(w)||/||yM(w)||\leq\]
\[
2\sum_{k\in S}|x_k-y_k|
\sum_{j\in S} \sum_{w\in {\cal W}} (M(w) )_{k,j} =
 2||x-y||\]
which combined with (\ref{G3eq11}) implies
that 
\[|T_{{\cal M}}u(x)-T_{{\cal M}}u(y)| \leq \] 
\[
||u||\cdot||x-y||+2\gamma(u)||x-y||
\]
from which (\ref{C5.eq0}) follows.

We now prove (\ref{R3eq1}). Thus let \(u \in Lip[K]\). That
(\ref{R3eq1}) holds if \(u\) is a constant is trivially true. Thus
assume that \(u \in Lip[K]\) is not a constant. Then \(\gamma(u) \not = 0\).
Now to prove (\ref{R3eq1}) set \(v=u/\gamma(u)\), set \(osc(v)= 
\sup \{|v(x)-v(y)|: x,y \in K\}\) and define \(v_0=v - osc(v)/2 - 
\inf\{v(x):x \in K\}\). Clearly \(\gamma(v)=\gamma(v_0)=1\).
 Since \(\sup\{||x-y||: x, y \in K\} =2\) it is also clear that
\(||v_0|| \leq 1\). 
Hence, by (\ref{C5.eq0}) follows 
\[
\gamma(T_{{\cal M}}v_0)\leq 3
\]
and hence
\[
\gamma(T_{{\cal M}}u)=\gamma(T_{{\cal M}}v_0)\gamma(u)\leq 3 \gamma(u)
\]
and thereby (\ref{R3eq1}) is proved. That (\ref{G3eq0}) also holds
follows trivially from (\ref{R3eq1}).
\(\;\;\Box\)

Since the right hand side of 
(\ref{R3eq1}) is independent of \(S, {\cal M}\) and \(P\) 
the following corollary follows immediately from 
Corollary \ref{cor3.4}.
\begin{corollary}\label{cor3.6} 
Let \(S\) be a denumerable set,
 let \({\cal M} \in {\cal G}(S)\) and
let \(u \in Lip[K]\).
Then
\begin{equation}\label{R3eq5a}
 \gamma(T^n_{{\cal M}}u) \leq 3 \gamma(u),\;\; n=1,2,... \;\;.
\end{equation}
\end {corollary}
{\em Proof}. Let \(u \in Lip[K]\). 
From 
Corollary \ref{cor3.4}
follows that 
\(\gamma(T_{{\cal M}}^nu)=
\gamma(T_{{\cal M}^n}u)\) and then (\ref{R3eq5a}) follows 
from Lemma \ref{lem3.5}. 
\(\;\;\Box
\)

\section{ A few facts about the Kantorovich metric.}
Let \({\cal Q}(K)\) denote the set 
of nonnegative, finite, Borel measures 
on \((K,{\cal E})\) with positive total mass. 
For \(\mu \in {\cal Q}(K)\) 
we write \(||\mu|| = \mu(K) \). For \(r>0\) we define 
\({\cal Q}_r(K)=\{\mu \in {\cal Q}(K): ||\mu|| = r \}\).
If both \(\mu,\nu \in {\cal Q}_r(K)\) 
we define, for any \(r >0,\)  
\begin{equation}\label{R3eq16c}
d_{K} (\mu, \nu) =
\sup\{ \int_{K} u(y)\mu(dy)-\int_{K}u(y)\nu(dy):
\; u \in Lip_1[K]\;\;\}.
\end{equation}
We call 
\(d_{K} (\mu, \nu)\) the Kantorovich distance between
\(\mu\) and \(\nu\).
Note that if \(\mu,\nu \in {\cal Q}_r(K)\) then 
both \(\mu/r\) and  \(\nu/r\) belong to \({\cal P}(K)\) and 
\begin{equation}\label{R3eq16d}
d_K(\mu, \nu)    =
rd_K(\mu/r,\nu/r).
\end{equation}
Note also that 
\[d_{K} (\mu, \nu) =
\sup\{ \int_{K} u(y)\mu(dy)-\int_{K}u(y)\nu(dy):
\; u \in Lip_1[K], ||u||\leq 1\;\;\}
\]
since \(\sup\{||x-y||: x,y \in K\}= 2.\)

That \(d_K(\cdot,\cdot)\) determines a metric on \({\cal P}(K)\) is
well-known, 
(see e.g. \cite{Dud02}, Chapter 11, section 3)
and from (\ref{R3eq16d}) follows that \(d_K\) determines a metric
also on \({\cal Q}_r(K)\) for any \(r > 0\).
We shall call the metric \(d_{K}(\cdot,\cdot)\) the 
{\em Kantorovich metric}.

For \(x \in K\) we let \(\delta_x\) denote the
probability measure in \({\cal P}(K) \) such that 
\(\delta_x(\{x\})=1.\)
From the definition of \(d_K(\cdot,\cdot)\) it readily follows that 
\[
d_K(\delta_x,\delta_y) =\delta(x,
y) = ||x-y||.
\]

If \({\cal P'} \subset {\cal P}(K) \) and \(\mu \in {\cal P}(K)\)
we define 
\[
d_K(\mu,{\cal P'}) =\inf \{d_K(\mu,\nu) : \nu \in {\cal P'}\}.
\]

It is well-known that the Kantorovich metric 
\(d_{K}\) on \({\cal Q}_r(K)\)
can be defined in another way. 
Let 
\(K^2= K\times K\), let
\({\cal E}^2= {\cal E}\otimes {\cal E}\),
let \(r >0\) and let \({\cal Q}_r(K\times K)\) denote the
set of nonnegative measures on 
\((K^2,{\cal E}^2)\)
 with total mass equal to \(r\).
For any two measures  \(\mu, \nu \in {\cal Q}_r(K)\)
we let
\(
{\tilde{\cal Q}}_r(\mu,\nu)\) denote the subset of 
\({\cal Q}_r(K\times K)\) consisting of those 
measures \(\tilde{\mu}(dx,dy)\)
such that 
\[\tilde{\mu}(A,K)= \mu(A), \;\;\;\;\;\;\;
\forall A \in {\cal E},\]
and
\[\tilde{\mu}(K,B)= \nu(B), \;\; \;\;\;\;\;\;\;\;\;
\forall B \in {\cal E}.\]
Then
\begin{equation}\label{R3eq24} 
d_{K}(\mu,\nu) =\inf \{ \int_{K\times K} \delta(x,y)
\tilde{\mu}(dx,dy): \; \tilde{\mu}(dx,dy)\in {\tilde{\cal Q}}_r(
\mu,\nu) \}.\end{equation}

A proof of the equality between (\ref{R3eq16c}) and (\ref{R3eq24})
when \(r=1\) can be found in 
\cite{Dud02}, section 11.8, and the equality between
(\ref{R3eq16c}) and (\ref{R3eq24}) for \(r\neq 1\) follows then by
using the relation (\ref{R3eq16d}).
The proof of the fact that the two definitions of \(d_K\) give the
same value 
goes back to
L.V. Kantorovich (see \cite{Kan42}).
For a short overview of the
Kantorovich metric and some applications,
see \cite{Ver06}.

Having introduced the metric \(d_K\) the following
corollary to Corollary \ref{cor3.6}
follows  immediately.

\begin{corollary} \label{cor4.1}
Let \(S\) be a denumerable set, 
let 
\({\cal M}\in G(S)\) 
and let \(\mu\) and \(\nu\) be two arbitrary measures in 
\({\cal P}(K)\). 
Then 
\[
d_{K}({\breve P}_{{\cal M}}^n\;\mu,
{\breve P}_{{\cal M}}^n\;\nu)
\leq 3 d_{K}(\mu,
\nu)\;\;for\;\; n=1,2,...\; .
\]
\end{corollary}
\begin{rem}
It is not difficult to construct an example that shows that
the constant 3, 
can not be replaced by a constant strictly less than 2.
\end{rem}

We end this section stating yet another 
lemma which is a simple consequence
of  well-known results from the general theory
on probability measures on complete, separable, 
metric spaces. See e.g \cite{Dud02} or \cite{Par67}, Chapter 2,
  section 6. 
That \(K\) is a complete, separable, metric space
follows because it is  a closed subset of \(l_1(S)\).

\begin{lemma}\label{lem4.2}
For every \(\mu \in {\cal P}(K)\) and every \(\epsilon >
  0\)
we can find an integer \(N\), a sequence \(\{\alpha_k,\; k=1,2,...,N\}\)
of real positive numbers satisfying 
\(
\sum_{k=1}^N \alpha_k = 1,
\)
and a sequence 
\(\{x_k, k=1,2,...,N\}\) of elements in \(K\) such that
if we define 
\( \nu = \sum_{k=1}^{N} \alpha_k \delta_{x_k} \)
then
\(\;
d_{K} (\mu, \nu) < \epsilon.
\;\)
\end{lemma}

\section{ On probability measures with equal barycenter.}

We start with a lemma that gives a lower bound of the 
Kantorovich distance. 

If \(\mu \in {\cal Q}(K)\) we define 
\({\overline b}(\mu)\) in the same way as when \(\mu \in {\cal P}(K)\) 
that is by the formula 
(\ref{i1eq12}) (see subsection 1.6).
\begin{lemma}\label{lem5.1} 
Let \(S\) be denumerable, let \(r > 0\) and 
let \(\mu, \nu  \in {\cal Q}_r(K)\).
Then 
\[||{\overline b}(\mu)-{\overline b}(\nu)|| \leq 
d_K(\mu,\nu). 
\]
\end{lemma}
{\em Proof.} 
Set 
\[a={\overline b}(\mu), \;\;and \;\;b={\overline b}(\nu).
\]
For 
\(i \in S\) 
define 
\[ u_i(x)= (x)_i.\]
By definition
\[
\int_{K}u_i(x)\mu(dx)=(a)_i \;\;and \;\; \int_{K}u_i(x)\nu(dx)=(b)_i
\;\;\;for \;\;i \in S. \]
Let \(sgn(\cdot)\) denote the so called sign-function defined by
\[sgn(t)= 1,  \;if \; \; t > 0,\;\;\;\;sgn(0)=0,\;\;and 
\;\;\;sgn(t)= -1,  \;if \; \; t < 0\;.\]

For \(i \in S \) we define
\[
\epsilon_i =  sgn(\int_{K}(x)_{i}\mu(dx)-
\int_{K}(x)_{i}\nu(dx)),  
\] 
and define \(v = \sum_{i \in S} \epsilon_i u_{i}\).
Obviously 
\[|v(x)-v(y)|\leq \sum_{i\in S}|(x)_i-(y)_i| =||x-y||
\]
and therefore \(\gamma(v) \leq 1\).
Furthermore 
\[
\int_{K}v(x)\mu(dx)- 
\int_{K}v(x)\nu(dx) =
 \sum_{i \in S}
\epsilon_i( \int_{K}(x)_{i}\mu(dx)-
\int_{K}(x)_{i}\nu(dx) \;)
 = \]
\[\sum_{i \in S}
|(a)_{i}-(b)_{i}| =||a-b||
\]
which implies that \(d_K(\mu,\nu)\geq ||a-b||\) since
\(v \in Lip_1[K]\).
\(\;\Box\)

\begin{lemma}\label{lem5.2} 
Let \(S\) be a denumerable set, 
let \(N\) be a positive integer 
let \(\xi_k, \;k=1,2,...,N\) be vectors in \(K\),
let \(\beta_k>0, \; k=1,2,...,N\) and  define the measure 
\(\varphi \in {\cal Q}(K)\) by 
\[
\varphi = \sum_{k=1}^N\beta_k\delta_{\xi_k}.
\]
Define the vector \(a\) by
\( a = \sum_{k=1}^N\beta_k\xi_k,
\)
let \(b=((b)_i, i \in S)\) be a vector with nonnegative coordinates
and such that
\(||b|| =||a||.\)

Then there exist vectors \(\zeta_k,\;\; k=1,2,...,N,\;\) in \(\;K\)
such that 
\(b= \sum_{k=1}^N\beta_k\zeta_k,\)
and such that if we define 
\(
\Psi =  \sum_{k=1}^N \beta_k \delta_{\zeta_k},
\)
then
\[
d_{K}(\varphi, \Psi) = 
\sum_k \beta_k ||\xi_k - \zeta_k|| = ||a-b||.
\]
\end{lemma}
{\em Proof.} Let 
\(
\varphi = \sum_{k=1}^N\beta_k\delta_{\xi_k}
\) be the given measure. If \(\zeta_k \in K, k=1,2,...,N\)
and 
\(
\;\Psi = \sum_{k=1}^N\beta_k\delta_{\zeta_k},\;
\)
then an upper bound for the Kantorovich distance
\(d_K(\varphi,\Psi)\) is obtained simply by the
estimate 
\(
\sum_{k=1}^N\beta_k||\xi_k-\zeta_k||
\)
since if we define the measure \({\tilde \varphi}\)
on \((K^2, {\cal E}^2)\) by
\[
{\tilde \varphi}(\{(\xi_k,\zeta_k)\}) = \beta_k, \;k=1,2,...,N
\]
then 
\[\tilde{\varphi}(A,K)= \varphi(A), \;\;
\forall A \in {\cal E},\;\;
and
\;\;
\tilde{\varphi}(K,B)= \Psi(B), \;\; 
\forall B \in {\cal E}\]
and therefore by (\ref{R3eq24}) follows that
\[
d_{K}(\varphi,\Psi) \leq
 \int_{K\times K} \delta(x,y)
\tilde{\varphi}(dx,dy) =
\sum_{k=1}^N\beta_k||\xi_k-\zeta_k||
\]
where we have used the fact that \(d_K(\delta_x,\delta_y)=||x-y||\).

From this observation and Lemma \ref{lem5.1} it follows that what we have to
do is to show that if \(b\) is a  given vector with 
nonnegative coordinates such that \(||b|| =
||a||\)
then we can find vectors 
\(\zeta_k,\;\; k=1,2,...,N,\;\) in \(\;K\) such that
\newline
1)
\[
\sum_{k=1}^N \beta_k\zeta_k = b
\]
2)
\[
\sum_{k=1}^N \beta_k||\xi_k-\zeta_k|| = ||a-b||.
\]
 
That this is possible when \(N=1\)  
is easily proved. Simply define
\[\zeta_1 = b/\beta_1.\]

The case when \(a=b\) is trivial. Just take 
\(\zeta_k=\xi_k, k=1,2,...,N.\)
In the remaining part of the proof
we therefore assume that \(a \neq b\). 

Let us now assume 
that we have 
proved the conclusion of the lemma
for \(N=M-1\) and let us prove the conclusion for \(N=M\), where 
\(M \geq 2\).

We shall have use of the following simple lemma which we state 
without proof.
\begin{lemma} \label{lem5.3} 
Let \(x,y,z \in l_1(S)\) have nonnegative coordinates
and assume that \(x+y=z\). Then \(||x||+||y||=||z||\).
\end{lemma}
 
Next let us define the sets \(S_1\), \(S_2\) and 
\(S_3\) by
\[
S_1 = \{ i \in S : (a)_i > (b)_i\;\},
\]
\[
S_2 = \{ i \in S: (a)_i < (b)_i\;\},
\]
and
\[
S_3 = \{ i \in S : (a)_i = (b)_i\;\}.
\]
Since \(a \neq b\) and \(||a||=||b||\)
 both the sets \(S_1\) and \(S_2\) are nonempty.
Let us set 
\[
\Delta=\sum_{i \in S_1} ((a)_i-(b)_i).
\]
Then clearly 
\[
\Delta= \sum_{j \in S_2} ((b)_j-(a)_j) = ||a-b||/2.
\]

Next let us consider the vector \(\xi_M\). We define
\[
R_1 = \{ i \in S_1 : (\xi_M)_i > 0 \}.
\]

Assume first that \(R_1= \emptyset\). Then define
the vector \(a_1\), the vector \(\zeta_m\) and the vector \(b_1\), by
\[a_1= a - \beta_M \xi_M,\;\;
\zeta_M = \xi_M,
\;\;
and
\;\;
 b_1 = b - \beta_M \zeta_M.\]
Since \(R_1=\emptyset\) it follows easily 
that the vector \(b_1\) has nonnegative
coordinates. Therefore by Lemma \ref{lem5.3} follows that
\(||a||=||a_1||+\beta_M||\xi_M||\) and also 
\(||b||=||b_1||+\beta_M||\zeta_M||\). Since 
\(\zeta_M = \xi_M\) and \(||a||=||b||\) it follows that 
\[||a_1||=||b_1||\;\; and\;\; 
||a-b||=||a_1-b_1||+ \beta_M||\xi_M- \zeta_M||.\]
 The conclusion of
Lemma \ref{lem5.2} now follows by the induction hypothesis.

It remains to consider the case when the set \(R_1\) is
non-empty.
Roughly speaking what we shall do is to move as large part
as possible from a coordinate \((\xi_M)_i\), \(i \in R_1\),
to a coordinate 
\((\xi_M)_j, j \in S_2\).

Hence, we claim that there exist 
nonnegative numbers \(t_{i,j}, \; i \in R_1, \;j \in S_2\)
with the following properties:
\[
\sum_{j \in S_2}t_{i,j} = 
\min \{ \beta_M (\xi_M)_i,((a)_i-(b)_i)\}, \;\;\;\; \forall i \in R_1
\]
and
\[
\sum_{i \in R_1} t_{i,j} \leq ((b)_j -(a)_j), \;\;\; \forall j \in S_2.
\]
That such a set \(\{ t_{i,j} : i \in R_1, \;\;j \in S_2\}\) exists
follows from the following two observations:
\newline
1) 
\[ \sum_{j \in S_2} ((b)_j-(a)_j) = \Delta , 
\]
2) 
\[\sum_{i \in R_1} \min\{\beta_M (\xi_M)_i,((a)_i-(b)_i)\} \leq
\sum_{i \in S_1} ((a)_i-(b)_i)\}  =
\Delta.
\]

We now simply define the vector \(\zeta_M\) by

\begin{equation}\label{G4eq1}
(\zeta_M)_{i} = (\xi_M)_i 
- (\beta_M)^{-1}\sum_{j \in S_2}t_{i,j}, \;\; i \in R_1
\end{equation}

\begin{equation}\label{G4eq2}
(\zeta_M)_j = (\xi_M)_j + 
(\beta_M)^{-1}\sum_{i \in R_1}t_{i,j}, \;\; j \in S_2
\end{equation}
\begin{equation}\label{G4eq3}
(\zeta_M)_i = (\xi_M)_i , \;\;if \;\;i \not \in R_1 \cup S_2.
\end{equation}
Since by definition
\[
\sum_{j \in S_2}t_{i,j} \leq \beta_M(\xi_M)_i, \;\;\;\forall i \in R_1
\]
it is clear that \((\zeta_M)_i\geq 0,\; \forall i \in R_1.\)
That also \((\zeta_M)_i \geq 0 \) if \( i \not \in R_1\) follows
from the fact that \((\xi_M)_i \geq 0\) for all \(i \in S\).
Hence 
\[
(\zeta_M)_i \geq 0, \;\forall i \in S
\]
and that \(||\zeta_M|| =||\xi_M|| =1\) now follows easily from
(\ref{G4eq1}), (\ref{G4eq2}) and (\ref{G4eq3}).

Now define 
\(a_1= a - \beta_M \xi_M \;\;and \;\; 
b_1 = b - \beta_M \zeta_M. \)
Obviously \((a_1)_i\geq 0,\;\forall i \in S\;\) 
and that also \((b_1)_i \geq 0\) for all \(i \in
S\)
is easily proved. 
First suppose \(\;i \in S_1\setminus R_1 \).
Then \((\xi_M)_i = 0\) and hence
\(
(b_1)_i = (b)_i \geq 0.
\)
Next suppose that \(\;i \in  R_1 \).
Then 
\[
(b_1)_i = (b)_i - \beta_M (\xi_M)_i + \sum_{j \in S_2}t_{i,j}
\geq \min \{ (b)_i, (a_1)_i\} \geq 0.
\]
If \(i \in S_3\) then 
\(
(b_1)_i=(a_1)_i \geq 0,  \;\)
and 
finally
if \(j \in S_2\) 
then 
\[
(b_1)_j = (b)_j - \beta_M (\zeta_M)_j \geq
 (a)_j-\beta_M (\xi_M)_j =(a_1)_j \geq 0.
\]

Hence \(||a_1||=||a||-\beta_M||\xi_M|| = 
||b||-\beta_M||\zeta_M|| =||b_1|| \) 
because of Lemma 5.3.
Furthermore from the definition of \(\zeta_M\) it follows easily that
\[
||a_1-b_1|| = 
2 \sum_{j \in S_2}((b)_j-(a)_j) -
2\sum_{i \in R_1, \;j \in S_2}t_{i,j} \]
and since 
\(
||a-b||=2 \sum_{j \in S_2}((b)_j-(a)_j) 
\)
and 
\(
2\sum_{i \in R_1, \;j \in S_2}t_{i,j} =
\beta_M ||\xi_M-\zeta_M||
\)
because of (\ref{G4eq1}) and (\ref{G4eq2}) we conclude that
\[
||a_1-b_1|| = ||a-b||-\beta_m ||\xi_M-\zeta_M||.
\]
Therefore,
just as was the case when \(R_1=\emptyset\), the
conclusion of Lemma 5.2 now follows by the induction hypothesis. 
\(\;\;\Box\)

The following two results are simple consequences of
Lemma \ref{lem5.1},
Lemma \ref{lem5.2}, Lemma \ref{lem4.2} and the triangle inequality. 
\newline
\begin{rem}
When \(S\) is finite and consequently \(K\) is compact,
then the conclusion of Corollary \ref{cor5.4} below
is well-known from the general theory on barycenters. (See e.g 
Proposition 26.4 in \cite{Cho69}.)
\end{rem}
\begin{corollary}\label{cor5.4} 
Let \(q \in K.\)
For every \(\mu \in {\cal P}(K|q)\) and every \(\epsilon >
  0\)
we can find an integer \(N\), a sequence \(\{\alpha_k, k=1,2,...,N\}\)
of real positive numbers satisfying 
\(
\sum_{k=1}^N \alpha_k = 1,
\)
and a sequence 
\(\{x_k, k=1,2,...,N\}\) of elements in \(K\) such that
if we define 
\( \nu = \sum_{k=1}^{N} \alpha_k \delta_{x_k} \)
then
\(
\nu \in {\cal P}(K|q)
\)
and 
\[
d_{K} (\mu, \nu) < \epsilon.
\]
\end{corollary}
\begin{corollary} \label{cor5.5}
Let \(\mu \in {\cal P}\) and  let \(q \in K\).
Then
\[
d_K(\mu, {\cal P}(K|q)) =||{\overline b}(\mu) - q || .
\]
\end{corollary} 
From Corollary \ref{cor5.4} and Corollary \ref{cor5.5}  also follows:
\begin{theorem}\label{thm5.6}
Let \(q \in K\). Then \({\cal P}(K|q)\) is a tight set.
\end{theorem}
{\em Proof.} It is well-known that \({\cal P}(K|q)\) is a compact set 
in the topology induced by the Kantorovich metric 
when the vector \(q\) has only finitely many nonzero coordinates.
From Corollary \ref{cor5.4} and Corollary \ref{cor5.5} 
it then easily follows that \({\cal P}(K|q)\) is a compact set in
the topology induced by the Kantorovich metric, also for arbitrary 
\(q\in K\). 
Therefore by \cite{Dud02}, Theorem 11.5.4 it  
follows that \({\cal P}(K|q)\) is tight
since  \((K,{\cal E})\) is a complete, separable,
metric space.
 \(\;\Box\)

We end this section with the following lemma to be used later.
For  \(i \in S\) and \(0 < \eta \leq 1\), we define
the set \(E_{i}(\eta)\) by
\begin{equation}\label{T4eq3}
E_{i}(\eta) =\{ x \in K : (x)_i \geq \eta \}.
\end{equation}

\begin{lemma}\label{lem5.7} 
Let \(S\) be denumerable, let \(i \in S,\)
let \(q \in K\) and suppose 
also that \((q)_i > 0\).
Then we can find a compact set \(C\)
such that for all \(\mu \in {\cal P}(K|q) \) 
\begin{equation}\label{S4eq20}
\mu(E_i((q)_i/2)\,\cap\, C) \geq (q)_i/3. 
\end{equation}
\end{lemma}
{\em Proof.} Let \(\mu \in {\cal P}(K|q) \). 
Since \(\int_{K}(y)_i\mu(dy) = (q)_i \) and \(0 \leq (y)_i \leq 1 \)
if \(y \in K\) one easily obtains the estimate
\[
\mu(E_i((q)_i/2)) \geq (q)_i/2
\]
for all \(\mu \in {\cal P}(K|q) \).
That we can find a compact set \(C\) such that 
(\ref{S4eq20}) holds, follows from the tightness of 
the set \({\cal P}(K|q) \). 
\(\;\Box
\)

\section{The barycenters of Markov chains induced by partitions.}

In this section we state some more results concerning 
barycenters. 

\begin{lemma}\label{lem6.1}
Let \(S\) be a denumerable set, let \(P\in PM(S)\) 
and let \({\cal M} \in {\cal G}(S,P)\) 
Then, for all \(x \in K\),
\begin{equation}\label{R5eq1}
{\overline b}({\breve P}^n_{{\cal M}}\delta_x) = xP^n, \;n=1,2,....
\end{equation}
\end{lemma} 
{\em Proof.} For an arbitrary \(i \in S\) 
define \(u_i\in C[K]\) by
\(
u_i(x)=(x)_i.
\)
It then follows that 
\[
\langle u_i, {\breve P}_{{\cal M}}\delta_x \rangle =
T_{{\cal M}}u_i(x)=\]
\[
\sum_{w\in {\cal W}_{{\cal M}(x)}}((xM(w))_i/||xM(w)||)
\cdot||xM(w) || =
\sum_{w\in {\cal W}_{{\cal M}(x)}}(xM(w))_i =(xP)_i\]
from which follows that (\ref{R5eq1}) holds for \(n=1\). 
That (\ref{R5eq1}) also holds for \(n \geq 2\) now follows from the
fact that \({\breve P}^n_{{\cal M}} = {\breve P}_{{\cal M}^n} \)
(see Corollary \ref{cor3.4}) and the fact that \({{\cal M}^n} \) 
is a partition of \(P^n\) (see (\ref{eq1.4.12})).
\(\;\;\Box\)

The following 
result is also easily proved if one uses Corollary \ref{cor5.4}.
\begin{theorem} \label{thm6.2}
Let \(S\) be a denumerable set, let \(P\in PM(S)\),  
 let \({\cal M} \in {\cal G}(S,P)\),
let \(\pi \in K\) and suppose that
\(
\pi =\pi P.
\)
Then, 
\[
{\breve P}_{{\cal M}}\mu \in {\cal P}(K|\pi), 
\;\;\forall \;\mu \in {\cal P}(K|\pi).
\]
\end{theorem}
{\em Proof.} 
First assume that \(\mu \in {\cal P}(K|\pi)\) can be written
\begin{equation}\label{R5eq2}
\mu = \sum_{k=1}^N \alpha_k \delta_{y_k}.
\end{equation}
Now, let \( i \in S\) be chosen arbitrarily, and let 
the function \(u \in C[K]\) be defined by \(u(x)=(x)_i\).
We obtain
\[
\langle u, {\breve P}_{{\cal M}}\mu \rangle =
\langle T_{\cal M}u,\mu \rangle =
\sum_{k=1}^N \alpha_k T_{\cal M}u(y_k) =
\]
\[
\sum_{k=1}^N \alpha_k \sum_{w\in {\cal W}_{{\cal M}}(y_k)}
 u(y_kM(w))/||y_kM(w)||)
\cdot ||y_kM(w)|| = \]
\[
\sum_{k=1}^N \alpha_k \sum_{w\in {\cal W}} (y_kM(w))_i 
=
(\pi P)_i = (\pi)_i
\]
and thereby we have proved the assertion of the theorem when 
the measure \(\mu\) can be written as in (\ref{R5eq2}).

That the conclusion is true for arbitrary \(\mu \in {\cal P}(K|\pi)\)
is easy to prove if one uses Corollary \ref{cor5.4}.
\( \;\Box \)

From Theorem \ref{thm6.2} and Theorem \ref{thm5.6}  
we now immediately obtain the following 
tightness result. 

\begin{theorem}\label{thm6.3}  Let \(S\) be a denumerable set, let
 \(P\in PM(S)\), let 
\({\cal M} \in {\cal G}(S,P)\), let  
\(\pi \in K\) and suppose that
\(\pi P= \pi\). 
Then, for all \(\;\mu \in {\cal P}(K|\pi)\)
\[
\{{\breve P}_{{\cal M}}^n \mu, \;\;n=1,2,...\}
\] 
is a tight sequence.
\end{theorem}

\section{ An auxiliary theorem for Markov chains in 
complete, separable, metric spaces.}

In this section we shall state a limit theorem for
Markov chains  in a complete, separable, metric space, which we
in the next section shall apply to  Markov chains generated
by tr.pr.fs induced by partitions of tr.pr.ms. 

In this section 
\((K,{\cal E})\) will denote an arbitrary complete,
separable, metric space, with metric \(\delta\) and 
\(\sigma-algebra\)  \({\cal E}\). 
Other notations will be the same as 
before.

Let  \(Q: K\times {\cal E}\rightarrow [0,1]\) 
be a tr.pr.f
on  \((K, {\cal E})\), and let 
\({ Q^n}: K\times {\cal E} \rightarrow [0,1], \;n=1,2,...\) 
be the sequence of a 
tr.pr.fs
defined
recursively
by
\[ Q^1(x, B)= Q(x,B), \;\;\;\;x\in K,\;\; B \in {\cal E}\]
\[ Q^{n+1}(x,B)=\int_{K}Q^n(y,B)\,Q(x,dy), \;
\;x\in K,\;\;B\ \in {\cal E}.\]
We let \( T: B[K] \rightarrow B[K] \)  denote the transition  operator 
associated  to \(Q\)
defined as usual by 
\(
Tu(x) = \int_{K}u(y)Q(x,dy).
\)
We define \(T^0u(x)=u(x)\). Note that 
\[
osc(T^{n+1}u) \leq 
osc(T^{n}u), \; \; \; n=0,1,2,...,  \; \; \; u\in B[K].
\]

When stating  and proving the forthcoming theorem we shall use the
notion {\em shrinking property}. (See (\ref{eq2.15}), for
  the 
definition of \(osc_A(u)\).)

\begin{definition} \label{def7.1}
Let \(Q\) be a tr.pr.f, and let \(T\) be the
associated transition operator.
If 
for every \(\rho > 0\) there exists a number 
\(\alpha \), \(\;0 < \alpha< 1\)  such that 
for every nonempty, compact set \(A \subset K\),
every \(\eta > 0\) and every \(\kappa > 0\), 
there exist an integer \(N\) and another 
nonempty, compact set \(B \subset K\) such that, 
if the integer \(n\geq N\) 
then for all \(u \in Lip[K]\)
\[
osc_A(T^{n}u)\leq
\eta \gamma(u) + \kappa osc(u)+ 
\alpha \rho \gamma(u) + 
(1-\alpha)
osc_B(T^{n-N}u)
\]
then we say that \(Q\) has the  {\bf shrinking property}.
We call \(\alpha \) a
{\bf shrinking number} associated to \(\rho\).
\end{definition}

We now first state and prove the following lemma.
\begin{lemma} \label{lem7.2} 
Suppose that the tr.p.f \(\;Q\) 
has the shrinking property.
Then for every nonempty compact set \(C \subset K\)
and every   
\(u \in Lip[K]\) 
\begin{equation}\label{S6eq0}
\lim_{n \rightarrow \infty} 
\sup_{x,y \in C}|\int_{K}u(z) Q^n(x,dz) - \int_{K}u(z) Q^n(y,dz)| = 0.
\end{equation}
\end{lemma}
{\em Proof.} 
Let \(C\) be a given nonempty, compact set. Set \(A_0=C.\)
Let also \(\epsilon > 0\) be
given. In order to prove the lemma it suffices to prove that, for
every
\(u \in Lip[K]\), we can find an integer \(N\) such that 
\[
osc_{A_0}(T^nu)
< 4\epsilon 
\]
if \(n \geq N\).

Thus let \(u \in Lip[K]\) be given.
Obviously (\ref{S6eq0}) holds if \(\gamma(u)=0\). Therefore we may
assume that \(\gamma(u) >0\), which also implies that \(osc(u)>0\).  

We now define
\begin{equation}\label{eq7.51}
\rho = \epsilon/\gamma(u).
\end{equation}
Next let \(\alpha > 0\) be a shrinking  number associated to
\(\rho\). Define the integer \(M\) by
\[
M=\min \{ m : (1-\alpha)^{m} < \epsilon/osc(u) \}.
\]
We now choose \(\eta_k=(1/2)^k\epsilon/\gamma(u),\; k=1,2,...,M \) 
and
\(\kappa_k=(1/2)^k\epsilon/osc(u),\; k=1,2,...M\). 
From the shrinking property follows that we can find
a sequence of
non-empty compact sets \(A_k,\;k=1,2,...,M\) and a sequence
of integers \(N_k\;k=1,2,...,M\) such that if \(m\geq N_k\) then 
\[osc_{A_{k-1}}(T^mu) \leq \eta_k\gamma(u)
+\kappa_k osc(u) + 
\alpha \gamma(u) \rho + (1-\alpha) osc_{A_k}(T^{m-N_k}u).
\]
Defining \(N = \sum_{k=1}^M N_k\) it follows easily by iteration 
and the definition of \(\rho\) (see (\ref{eq7.51})) that
if \(n \geq N\) then
\[
osc_{A_0}(T^nu) \leq (1/2)\epsilon + (1/2)\epsilon +
\alpha\epsilon + (1-\alpha)osc_{A_1}(T^{n-N_1}u) \leq \]
\[(1/2+1/4)\epsilon + (1/2+1/4)\epsilon +
\epsilon\alpha(1+(1-\alpha))+
(1-\alpha)^2osc_{A_2}(T^{n-N_1-N_2}u) \leq ... \leq\]
\[
(\sum_{k=1}^M(1/2)^k)\epsilon +
(\sum_{k=1}^M(1/2)^k)\epsilon + 
\epsilon\alpha(\sum_{k=0}^{M-1}(1-\alpha)^k) +
(1-\alpha)^M osc_{A_N}(T^{n-N}u) < \]
\[
\epsilon + \epsilon + \epsilon + \epsilon
= 4\epsilon 
\] 
which was what we wanted to prove. \(\;\;\Box\)

\begin{theorem}\label{thm7.3}
Suppose that the tr.p.f \(Q\) is Feller continuous,
that \(Q\) has the 
shrinking property
and that there exists a point
\(x^* \in K\) such that the sequence 
\(
\{
Q^{n}(x^*,\cdot), n=1,2,...\}
\) 
is a tight sequence 
of probability measures. 
Then Q is asymptotically stable.
\end{theorem}
{\em Proof.} 
Firstly, using the fact that \(Q\) is Feller continuous 
and that 
\(
\{
Q^{n}(x^*,\cdot), n=1,2,...\}
\)
is a tight sequence it is well-known from the general theory on  Markov
 chains with complete separable metric state spaces
that there
exists at least one invariant measure for \(Q\). 

Secondly, using Lemma \ref{lem7.2} it is easy to prove - by contradiction -
that there is only one invariant measure, \(\nu\) say.
For suppose both 
\(\nu\) and \(\mu\) are invariant measures for \(Q\). Suppose \(||u||=1\),
and that \(a > 0 \) where \(a\) is defined by
\[ 
a= \int_{K}u(x)\mu(dx) - \int_{K}u(x)\nu(dx).
\]
Choosing the compact set \(C\) sufficiently large it is clear
that for 
\(n=0,1,2,... \) we
have 
\[ 
a= \int_{K}T^nu(x)\mu(dx) - \int_{K}T^nu(y)\nu(dy) <
 \]
\[
a/4 + \int_{C} \int_{C} |T^nu(x)-T^nu(y)|\mu(dx)\nu(dy) 
\]
and then from Lemma \ref{lem7.2} we can conclude that
\[
\int_{C} \int_{C} |T^nu(x)-T^nu(y)|\mu(dx)\nu(dy) < a/4
\] if \(n\) sufficiently large and hence \(a < a/2\) and we have
obtained our contradiction.

Thirdly, using Lemma \ref{lem7.2} again,  it is also easily proved that
\begin{equation}\label{S6eq51}
\lim _{n \rightarrow \infty}|T^{n+1}u(x^*)-T^nu(x^*)|=0 
\end{equation}
for \(u \in Lip[K]\).
Thus, let \(u \in Lip[K] \) be given and let also \(\epsilon > 0\) be given. 
Choose the compact set \(C\) so large that \((1-Q(x^*,C))osc(u) <
\epsilon/2\;\). Since 
\[ |T^{n+1}u(x^*)-T^nu(x^*)| =
|\int_{K} (T^{n}u(y)-T^nu(x^*))Q(x^*,dy)|\]
 it follows that
\[ |T^{n+1}u(x^*)-T^nu(x^*)| \leq
 |\int_{C} (T^{n}u(y)-T^nu(x^*))Q(x^*,dy)|+ \epsilon/2
\]
and since 
\[ |\int_{C} (T^{n}u(y)-T^nu(x^*))Q(x^*,dy)|< \epsilon/2\]
if \(n \) sufficiently large because of Lemma \ref{lem7.2}, statement 
(\ref{S6eq51}) follows.

From (\ref{S6eq51}) and the fact that 
\(
\{
Q^{n}(x^*,\cdot), n=1,2,...\}
\)
is a tight sequence it 
follows that 
\[
\lim _{n \rightarrow \infty}T^nu(x^*) = \int_K u(y)\nu(dy)
\] 
for all \(u \in Lip[K]\), and using  Lemma \ref{lem7.2} yet again
it 
also follows that 
\begin{equation}\label{S6eq5}
\lim _{n \rightarrow \infty}T^nu(x) = \int_K u(y)\nu(dy)
\end{equation}
for all \(u \in Lip[K]\) and all \(x \in K\).

Finally, that (\ref{S6eq5}) also holds for all \(u \in C[K]\) 
and all \(x \in K\) follows from the fact that
the set \(Lip[K]\) is measure determining in a separable, metric
space.
(See for example \cite{Dud02}, Theorem
11.3.3.)
\(\;\Box\)

We end this section introducing the following terminology.
\begin{definition} \label{def7.4}  
Let \(Q\) be a tr.pr.f and let \(T\) be the
associated transition operator.\newline
(i) If there exists a positive constant \(L\)
such that 
for every \(u \in Lip[K]\)
\[
|Tu(x)-Tu(y)| \leq L\gamma(u) \delta(x,y),\;\;\;  
\]
then we say that \(Q\) is {\bf Lipschitz continuous}. 
\newline
(ii)
If there exists a positive constant \(L\)
such that for
every \(u \in Lip[K]\)
\[
|T^nu(x)-T^nu(y)| \leq L\gamma(u) \delta(x,y),\;\;\;  n = 0,1,2,....  
\]
then we say that \(Q\) is {\bf Lipschitz equicontinuous}. 
\end{definition}

\begin{proposition}\label{prop7.5}
Let \(S\) be a denumerable set, let \({\cal M} \in {\cal G}(S)\) 
and let 
\(T_{{\cal M}}\) denote the transition operator on \(K\) induced 
by \({\cal M}\). Then \(T_{{\cal M}}\) is both Lipschitz continuous
and Lipschitz equicontinuous.
\end{proposition}
{\em Proof.} Follows immediately from Corollary \ref{cor3.6}.
\(\;\;\Box\)

\section{The proof of Theorem \ref{thm1.1}.}

We first repeat the formulation of Theorem \ref{thm1.1}.
\newline

{\bf Theorem \ref{thm1.1}} 
{\em Let 
\(S\) be a denumerable set, let \(P \in PM_{ae}(S),\)
let \(\pi \in K\) satisfy \(\pi P =\pi\),
let \({\cal M}\in G(S,P)\) 
and let \({\bf P}_{{\cal M}}\) be the tr.pr.f induced by  
\({\cal M}\).
Suppose also that 
Condition B holds. Then 
\({\bf P}_{{\cal M}}\) is asymptotically stable. 
}
\newline

{\em Proof.} 
Since \((K,{\cal E})\) is a complete, separable, metric space
it suffices to prove that \({\bf P}_{{\cal M}}\) satisfies 
the hypotheses of Theorem \ref{thm7.3}.
 
That \({\bf P}_{{\cal M}}\) is {\em Feller continuous} means that 
\begin{equation}\label{eq7.1}
u \in C[K] \Rightarrow T_{{\cal M}}u \in C[K].
\end{equation}
Since \(T_{{\cal M}}u(x) = \sum_{w \in 
{\cal W}_{{\cal M}(x)}}u(xM(w)/||xM(w)||)||xM(w)|| \) it it is clear 
that if the set \({\cal W}\) is finite then (\ref{eq7.1}) holds,
since for each \(w \in {\cal W}\) the function \(f_w:K 
\rightarrow {\tt R}\) defined by \(f_w(x)=||xM(w)||\) is continuous.
That  
(\ref{eq7.1}) also holds when \({\cal W}\) is an infinite set 
then follows by a simple truncation argument.
 
From 
Theorem \ref{thm6.3}
we furthermore conclude that 
\(
\{
{\bf P}_{{\cal M}}^{n}(\pi,\cdot), n=1,2,...\}
\) 
is {\em a tight sequence} 
since \({\cal M} \in 
{\cal G}(S,P)\), 
\(\pi P =\pi\) and  \(\delta_{\pi} \in {\cal P}(K|\pi)\).

It thus remains to show that 
\({\bf P}_{{\cal M}}\) 
has
the shrinking property. To simplify notations we shall throughout 
the rest of this
proof
denote the transition probability function 
\({\bf P}_{{\cal M}}\)  by \({\bf P}\), denote the
transition operator 
\(T_{{\cal M}}\) by \(T\) and the transition probability operator
\({\breve P}_{{\cal M}} \) by
   \({\breve P} \). 

Thus let \(\rho> 0 \) be given. What we have to do is to show
that we can find a number \(\alpha > 0\) such that for each
nonempty compact set \(A\) and each \(\eta > 0\) and each
\(\kappa > 0\)
we can find an integer \(N\) and another nonempty compact set \(B\)
such that for each \(u \in Lip[K]\)
\begin{equation}\label{S7eq10}
osc_A(T^nu) \leq
\eta \gamma(u) + \kappa osc(u)+ 
\alpha \rho \gamma(u) + 
(1-\alpha)osc_B(T^{n-N}u).
\end{equation}

Now, let also \(A\), \(\eta\) and 
\(\kappa\) be given where thus \(A\) is a nonempty compact set,
\(\eta > 0\) and \(\kappa > 0\).
Verifying the shrinking property will be
 done in three steps. In the first step we choose 
\(N_1\) so large that if \(n \geq N_1\) then the barycenter of 
\({\breve P}^n\delta_x \) is very close to \(\pi\) for every \(x \in
A\). The integer \(N_1\) will depend on \(A\) and \(\eta\).  
In the second step we use Condition B to determine 
a shrinking coefficient \(\alpha\) and an integer \(N_2\) -
only depending on \(\rho\) -  and use tightness of the set 
\({\cal P}(K|\pi)\) and the given number \(\kappa\) to determine
the compact set \(B\). In the third step we make the necessary
estimations in order to verify the shrinking property. 

{\bf Step 1.}
Since the set \(A\) given above is compact, it follows 
from (\ref{eq2.1.1}) 
that we can find an integer \(N_1\) such that for all \(z \in A\)
\[
||zP^n-\pi||< \eta/6
\]
if \(n \geq N_1\).
Now let \(x \in A \) and \(y \in A\) be given. Set 
\[ \mu_{n,x}(\cdot) =
{\bf P}^{n}(x,\cdot), \;\;
and
\;\;
 \mu_{n,y}(\cdot) =
{\bf P}^{n}(y,\cdot), \;\;n=1,2,...\; .
\]
From Lemma \ref{lem6.1} follows that 
\(
{\overline b}(\mu_{n,x}) = x P^n 
\)
and that 
\(
{\overline b}(\mu_{n,y}) = y P^n.
\)
Therefore, if \(n\geq N_1\), where \(N_1\) is defined as above,
we conclude that
\[
|| {\overline b}(\mu_{n,x})-\pi ||< \eta/6
\;\;and \;\;
|| {\overline b}(\mu_{n,y})-\pi ||< \eta/6.
\]
From Corollary \ref{cor5.5}  now follows that we can find two measures 
\(\nu_x\) and \(\nu_y\), both
in  \({\cal P}(K|\pi)\), such that if \(u \in Lip[K]\) then
\[
|\int_{K}u(z)\mu_{N_1,x}(dz) -  \int_{K}u(z)\nu_{x}(dz)| \leq \gamma(u)\eta/6
\]
and
\[
|\int_{K}u(z)\mu_{N_1,y}(dz) -  \int_{K}u(z)\nu_{y}(dz)|  
\leq \gamma(u)\eta/6.  
\]

From Corollary \ref{cor3.6}
(the Lipschitz equicontinuity property)
we also find that for \(m=0,1,2,...\) 
\[
|\int_{K}T^mu(z)\mu_{N_1,x}(dz) -  
\int_{K}T^mu(z)\nu_{x}(dz)| \leq 3 \gamma(u)\eta/6
 =\gamma(u)\eta/2     
\]
and similarly that
\[
|\int_{K}T^m_{{\cal M}}u(z)\mu_{N_1,y}(dz) -  
\int_{K}T^m_{{\cal M}}u(z)\nu_{y}(dz)| \leq 3 \gamma(u) \eta/6 =
\gamma(u)\eta/2.  
\]
Thus if \(n \geq N_1 \) we have
\[
|T^nu(x) - T^nu(y)|\leq 
\eta \gamma(u)+\]
\begin{equation}\label{eq7.2}
|\int_{K}T^{n-N_1}u(z)\nu_x(dz) - 
\int_{K}T^{n-N_1}u(z)\nu_y(dz)|.
\end{equation}
This concludes the first step.

{\bf Step 2.}
We shall next define a shrinking coefficient \(\alpha > 0\)
associated to the given number \(\rho\). To do this we shall use
Lemma \ref{lem5.7} and Condition B.

From Lemma \ref{lem5.7} it follows that for each \(i \in S\) we can find a
compact set
 \(C_i\) 
such that
for all
\(\mu \in {\cal P}(K|\pi)\)
\begin{equation}\label{S7eq80}
\mu(C_i\cap E_{i}((\pi)_{i}/2)) \geq (\pi)_{i}/3.
\end{equation}

From Condition B it follows that we can find  
an element \(i_0 \in S\), an integer \(N_0\)  and a sequence 
\(\{w_1,w_2,...w_{N_0}\}\) of elements in 
\({\cal W}\) depending on the set \(C_{i_0}\), such that if we set 
\(M(w_1)M(w_2)...M(w_{N_0})={\bf M}({\bf w^{N_0}})\) then 
\(
||e^{i_0}{\bf M}({\bf w^{N_0}}))||>0
\)
and for all \(x \in E_{i_0}((\pi)_{i_0}/2)\cap C_{i_0} \)
we  have
\begin{equation}\label{T7eq31}
||(x{\bf M}({\bf w^{N_0}})/||x{\bf M}({\bf w^{N_0}})|| -
e^{i_0}{\bf M}({\bf w^{N_0}})/e^{i_0}{\bf M}({\bf w^{N_0}})||)||< \rho/6.
\;\;\;
\end{equation}
Let us now define \(\alpha_1\) by
\begin{equation}\label{S7eq81}
\alpha_1=((\pi)_{i_0}/3)\cdot
((\pi)_{i_0}/2) \cdot ||e^{i_0}{\bf M}({\bf w^{N_0}})||,
\end{equation}
and let us define 
\[
\alpha =\alpha_1^2/2.
\]

Our aim is to verify (\ref{S7eq10}) with this choice of
\(\alpha\) and with \(N=N_1+N_0\).
In order to do this let us first set
\[
\nu^*_x = {\breve P}^{N_0}\nu_x,
\;\;\;
\nu^*_y = {\breve P}^{N_0}\nu_y,
\]
and let 
\({\tilde \nu}_{x,y}^*\) 
denote the
product measure on \((K^{2},{\cal E}^{2})\)
determined by 
\(\nu_x^*\) and
\(\nu_y^*\).

Furthermore let  us denote 
\[q_0 = e^{i_0}{\bf M}({\bf w}^{N_0})/
||e^{i_0}{\bf M}({\bf w}^{N_0})||\]
and
\[
D= \{z \in K :  \delta(z,q_0)<\rho/6\}. 
\]
Since \((z)_{i_0} \geq (\pi)_{i_0}/2 \) if 
\(z \in E_{i_0}((\pi)_{i_0}/2)\)
and 
\begin{equation}\label{T7eq91}
||xM|| \geq (x)_i||e^iM||, \;\; \forall i \in S, 
\end{equation}
if \(M\) is a nonnegative \(S\times S\) matrix and \(x \in K\), 
we conclude from (\ref{S7eq80}), (\ref{T7eq31}), (\ref{S7eq81}) 
and (\ref{T7eq91}) that
\[
\nu_x^*(D) \geq ((\pi)_{i_0}/3)\cdot ((\pi)_{i_0}/2)
\cdot ||e^{i_0}{\bf M}^{N_0}({\bf w}^{N_0})|| = \alpha_1
\]
and that the same inequality holds for \(\nu_y^*(D)\).

Since 
\(\nu_x^* = {\breve P}^{N_0}\nu_x\;\),  
\(\nu_y^* = {\breve P}^{N_0}\nu_y\) 
and
\(\nu_x, \nu_y \in {\cal P}(K|\pi)\),  
it follows from
Theorem \ref{thm6.2} that  \(\nu_x^*, \nu_y^* \in {\cal P}(K|\pi)\).
Since \({\cal P}(K|\pi)\) is a tight set 
it follows 
that we can find a compact set \(B\) independent of \(x,y \in A\)
and also of \(A\), 
 such that
\begin{equation}\label{eq7.4}
\nu_x^*(D\cap B) \geq \alpha_1/\sqrt 2,
\;\; and\;\;
\nu_y^*(D\cap B) \geq \alpha_1/\sqrt 2,
\end{equation}
and also that
\begin{equation}\label{eq7.5}
\nu_x^*(B) \geq 1 -\kappa/2
\;\;and\;\;
\nu_y^*(B) \geq 1 -\kappa/2.
\end{equation}
This concludes step 2.
\newline

{\bf Step 3.} Set \(N= N_1+N_0\). In this last step we shall estimate
\[|\int_{K}T^{n-N_1}u(z)\nu_x(dz) - 
\int_{K}T^{n-N_1}u(z)\nu_y(dz)|
\] when \(n \geq N\).
Set \(m=n-N\) and \(v= T^mu\). Then \(n-N_1= m+N_0.\) Hence
\[
|\int_{K}T^{n-N_1}u(z)\nu_x(dz) - 
\int_{K}T^{n-N_1}u(z)\nu_y(dz)| = \]
\[
|\int_{K}T^{m+N_0}u(z)\nu_x(dz) - 
\int_{K}T^{m+N_0}u(z)\nu_y(dz)|=
|\langle T^{m}u,\nu_x^*\rangle - 
\langle T^{m}u,\nu_y^*\rangle|=
\]
\[
|\int_{K^{(2)}} (v(z)-v(z')) {\tilde \nu}_{x,y}^*(dz,dz')|\leq
\]
\[|\int_{B_1} (v(z) - v(z')){\tilde\nu}_{x,y}^*(dz,dz')|+
|\int_{B_2} (v(z) - v(z')){\tilde\nu}_{x,y}^*(dz,dz')|+
\]
\[
|\int_{B_3} (v(z) - v(z')){\tilde\nu}_{x,y}^*(dz,dz')|
\]
where
\[
B_1= \{
(z,z') \in K^2: \delta(z,z')<\rho/3, z\in B, z' \in B\},
\]
\[B_2=\{(z,z')\in K_2 : \delta(z,z')\geq \rho/3, z\in B, z' \in B\},
\;and \; B_3= K^2\setminus (B_1\cup B_2). \]
From (\ref{eq7.5}) follows easily that
\begin{equation}\label{eq7.11}
|\int_{B_3} (v(z) - v(z')){\tilde\nu}_{x,y}^*(dz,dz')|
\leq osc(v)(1-\kappa).
\end{equation}
Furthermore
\[
|\int_{B_2} (v(z) - v(z')){\tilde\nu}_{x,y}^*(dz,dz')|
\leq osc_B(v){\tilde\nu}_{x,y}^*(B_2))\leq
osc_B(v)(1-{\tilde\nu}_{x,y}^*(B_1)).
\]
But if \(0 <a \leq b \leq 1 \), \(\epsilon > 0\) and 
\(\Theta > 0\), then by elementary calculations we find
\begin{equation}\label{eq7.6}
b \min \{\epsilon,\Theta\} + (1-b) \Theta \leq
a \epsilon + (1-a)\Theta.
\end{equation}
Since \( {\tilde\nu}_{x,y}^*(B_1^{(2)})\geq \alpha\) because
of (\ref{eq7.4}), 
it follows from (\ref{eq7.6}) that  
\[|\int_{B_1^{(2)}} (v(z) - v(z')){\tilde\nu}_{x,y}^*(dz,dz')|+
|\int_{B_2^{(2)}} (v(z) - v(z')){\tilde\nu}_{x,y}^*(dz,dz')|\leq\]
\[
\min \{\gamma(v)\rho/3,osc_B(v) \}\cdot {\tilde\nu}_{x,y}^*(B_1^{(2)}) + 
osc_B(v) (1-{\tilde\nu}_{x,y}^*(B_1^{(2)})) \leq \]
\[\alpha \gamma(v)\rho/3 +
(1-\alpha)osc_B(v) \]
which combined with (\ref{eq7.11})
implies that
\[
|\int_{K}T^{n-N_1}u(z)\nu_x(dz) - 
\int_{K}T^{n-N_1}u(z)\nu_y(dz)| \leq
\]
\begin{equation}\label{eq7.12}
osc(u)\kappa +
\alpha \gamma(u) \rho + (1-\alpha) osc_B(T^{n-N}u)
\end{equation}
where we also used the fact that \(\gamma(v)\leq 3\gamma(u)\)
because of Corollary \ref{cor3.6}.
That (\ref{S7eq10}) holds now follows by combining
 (\ref{eq7.12}) and (\ref{eq7.2}), and thereby Theorem \ref{thm1.1} 
is proved. 
\(\;\Box\)

\section{On Condition B.}

In this section we shall prove that Condition B1 (see
subsection 1.5) implies Condition B, we shall introduce
a notion we call
{\em localization}.
and present a theorem based on  Condition A. 

\begin{proposition} \label{prop9.1}
Let \(S\) be a denumerable set, let \(P \in PM_{ae}(S)\),
let \({\cal M}=\{M(w) \,: \, w \in {\cal W}\}\in {\cal G}(S,P)\),
and let \(\pi \in K\) satisfy
\(\pi P =\pi\).
Then,  if Condition B1 is satisfied, it follows that 
Condition B is also satisfied.
\end{proposition}
{\em Proof.}
In order to prove Proposition 9.1  we shall need the
following lemma in which we state an inequality
for matrices approaching a matrix in the set 
\({\tt W} \). (For the definition of the set 
\({\tt W} \) see subsection 1.6.) 
\begin{lemma}\label{lem9.2} 
Let \(u \in {\cal U}, \; v\in K\), and define
\(W=u^cv\). Let \(i_0\) be such that
\[
(u)_{i_0}> 0.
\]
Let \(\{ W_n, n = 1,2,... \}\) 
be a sequence of 
matrices of the same format as \(W\), and assume that 
\newline
1) for \(n=1,2,... \)
\[||W_n|| =1,\]
2) 
\begin{equation}\label{eq9.9}
\lim_{n \rightarrow \infty} ||e^i(W_n-W)|| = 0, \; i \in S.
\end{equation}
Then to every 
\(\eta, \; 0 < \eta < 1\)
and every nonempty compact set \(C\subset K\)
and  every \(\gamma > 0\),
there exists an integer  \( N=N_{\gamma,C,\eta}\)
such that if
\(x\in C\) and \(y \in C\) are such that
\( (x)_{i_0}\geq  \eta\) and
\( (y)_{i_0}\geq  \eta\)
then   \(||xW_n|| > 0\) and \(||yW_n|| > 0\) 
for all integers 
\(n \geq N\), and furthermore
\[
|| (xW_n/||xW_n|| - yW_n/||yW_n||) || < \gamma . 
\]
\end{lemma}
{\em Proof of Lemma \ref{lem9.2}.} Let \(C\) be a given nonempty compact set.
Since \(C\) is compact it readily follows from 
(\ref{eq9.9}) that to every \(\epsilon\)
we can find an integer \(N\) such that if \(n\geq N\) then
\begin{equation}\label{eq9.12}
|(||xW_n||-||xW||)|< \epsilon, \;\forall x \in C.
\end{equation}
Next let \(\eta\) be given and assume that \(x\in C\) and 
\((x)_{i_0}>\eta.\) 
Then again using (\ref{eq9.9}) it is easily proved
that if \(n\) is sufficiently large then
\begin{equation}\label{eq9.13}
||xW_n|| > (u)_{i_0}\eta/2, \;\;\;\forall x \in C
\;such\;\; that \;\; (x)_{i_0} \geq \eta.
\end{equation}
Finally, letting \(\gamma > 0\) also be given, by
 using the inequality (\ref{R3eq6}) together with
(\ref{eq9.12}), (\ref{eq9.13}) and the triangle inequality
we can conclude that
 \[
|| (xW_n/||xW_n|| - yW_n/||yW_n||) || < 
2 ||xW_n-yW||/
\max \{||xW_n||,||yW_n||\} \leq \]
\[
2||xW_n-yW_n||(2/(u)_{i_0}\eta) <
\gamma 
\]
if \(n\) is sufficiently large. \(\;\Box\)

We now continue the proof of  Proposition \ref{prop9.1}.
Let \(\rho > 0\) be given.
In order to prove Proposition \ref{prop9.1} we shall  prove  that 
we can find 
an element \(i_0 \in S\), such that if
\(C\) is a compact set such that
\begin{equation}\label{S5eq11}
\mu(C\cap E_{i_0}((\pi)_{i_0}/2))  \geq (\pi)_{i_0}/3
\end{equation}
for all \(\mu \in {\cal P}(K|\pi)\),
then we can find 
an integer \(N\) and an element 
\({\bf w^N} \in {\cal W}^N\)
such that 
\[
||e^{i_0}{\bf M}({\bf w^{N}})|| >0
\]
and 
\[
||(x{\bf M}({\bf w^{N}})/||x{\bf M}({\bf w^{N}})||
-e^{i_0}{\bf M}({\bf w^{N}})/||e^{i_0}{\bf M}({\bf w^{N}})||)||
< \rho,
\]
\begin{equation}\label{G9eq52}
\;\;\;\;\;\; \forall x \in E_{i_0}((\pi)_{i_0}/2)\cap C. 
\end{equation}
That there exists a compact set \(C\) such that 
(\ref{S5eq11}) holds 
follows from Lemma \ref{lem5.7}.

Since Condition B1 is satisfied
there exist a vector \(u \in {\cal U}\), a vector \(v \in K\), 
a sequence of integers \(\{n_1,n_2,...\}\), and a sequence
\(\{{\bf w_j^{n_j}}=\{w_{1,j},w_{2,j},...,w_{n_j,j}\}, \;j=1,2,...\}\)
of sequences 
%of elements in 
%\({\cal W}^{n_j}\) respectively,
such that 
\(||{\bf M}({\bf w_j^{n_j}})||>0,\; j =1,2,...\) 
and such that if we define \(W = u^cv\) then 
for all \(i \in S\) 
\begin{equation}\label{G9eq53}
\lim_{j\rightarrow \infty} 
||(e^{i}{\bf M}({\bf w_j^{n_j}})/
 ||{\bf M}({\bf w_j^{n_j}})||) - 
 e^{i}W|| = 0.
\end{equation}

Let us choose \(i_0 \in S\) such that 
\(
(u)_{i_0}>0\)
and let \(C\) be a compact set such that
(\ref{S5eq11}) holds for all \(\mu \in {\cal P}(K|\pi)\).
Since \((u)_{i_0}>0\) it follows that \(||e^{i_0}W||= 
(u)_{i_0}>0\). By (\ref{G9eq53}) then follows that 
\(||{\bf M}({\bf w_j^{n_j}})||>0,\; j =1,2,...\)
if we let the  enumeration \(n_1,n_2,...\) start with a sufficiently
large 
\(n_1\).
Since obviously
\[||({\bf M}({\bf w_j^{n_j}})/
 ||{\bf M}({\bf w_j^{n_j}})||)|| = 1 
\]
and \(C\cup\{e^{i_0}\}\) is a compact set,
it follows from Lemma \ref{lem9.2} that
if \(j\) is sufficiently large then

\[
||(x{\bf M}({\bf w^{n_j}})/||x{\bf M}({\bf w^{n_j}})||
-e^{i_0}{\bf M}({\bf w^{n_j}})/||e^{i_0}{\bf M}({\bf w^{n_j}})||)||
< \rho,
\]
\[
\;\; \forall x \in E_{i_0}((\pi)_{i_0}/2)\cap C
\]
and hence (\ref{G9eq52}) holds 
which was what we wanted to prove. \(\Box\)

Next, we call a  matrix  \(M\) {\bf subrectangular} if 
\[
(M)_{i_1,j_1}\neq 0 \;and  \;\;(M)_{i_2,j_2} \neq 0 \;\;
\Rightarrow 
\;\;
(M)_{i_1,j_2}\neq 0 \;\;and\;\;\; \;(M)_{i_2,j_1} \neq 0. 
\]
Before we state our next theorem we need one more notion.
As usual let \(S\) be a denumerable set, and let \({\cal M}=
\{M(w): w \in {\cal W}\} \in G(S).\)
We say that
\({\cal M}\) is {\bf localizing} if there exists a sequence
   \(w_1,w_2,...,w_n\) such that the number of non-zero columns
of the matrix \({\bf M}({\bf w^n})\) is {\em finite}. 

\begin{theorem} \label{thm9.3}  
Let 
\(S\) be a denumerable set, let \(P \in PM_{ae}(S)\),
let \({\cal M}=\{M(w): w \in {\cal W}\} \in G(S,P)\)
and let \({\bf P}_{{\cal M}}\) be the tr.pr.f induced by  
\({\cal M}\).
Suppose also that
\newline
1) \({\cal M}\) is localizing
\newline
2) Condition A is satisfied.

It then follows that 
\({\bf P}_{{\cal M}}\) is asymptotically stable.
\end{theorem}
{\em Proof.} From Proposition \ref{prop9.1} and Theorem \ref{thm1.1}
 follows that
it suffices to verify Condition B1.
In order to do this we follow closely the arguments used by
Kochman and Reeds
in their proof of Theorem 2 in \cite{KR06}. 

Since we have assumed that Condition A is satisfied
we can find  elements \(a_1,a_2,...,a_{N_1}\) in \({\cal W}\) such that
\({\bf M}({\bf a^{N_1}})\) is a non-zero subrectangular matrix.
Since \({\cal M}\) is localizing we can find elements
\(b_1,b_2,...,b_{N_2}\) in \({\cal W}\) such that
\({\bf M}({\bf b^{N_2}})\) has finitely many non-zero
columns. Choose \(i_1,j_1\) such that 
\(({\bf M}({\bf a^{N_1}}))_{i_1,j_1}>0\) 
and \(i_0,j_0\) such that 
\(({\bf M}({\bf b^{N_2}}))_{i_0,j_0}>0.\) 
Now, since \(P\) is irreducible there exist an integer \(N_3\) and 
 elements  \(c_1,c_2,...,c_{N_3}\) in \({\cal W}\) such that
%\({\bf c^{N_3}}\in {\cal W}^{N_3} \) such that
\(({\bf M}({\bf c^{N_3}}))_{j_1,i_0} > 0\) and also an integer
\(N_4\)
and elements
\(d_1,d_2,...,d_{N_4}\) in \({\cal W}\) such that 
%\({\bf d^{N_4}}\in {\cal W}^{N_4}\) such that 
\(
({\bf M}({\bf d^{N_4}}))_{j_0,i_1} >0.\)

Now let us define \(N=N_1+N_2+N_3+N_4\), define 
the element \({\bf w^{N}}=\{w_1,w_2,...,w_N\} \) by
\[
{\bf w^{N}} = ({\bf d^{N_4}},{\bf a^{N_1}},
{\bf c^{N_3}},{\bf b^{N_2}})\]
and define
\[
G = {\bf M}({\bf w^{N}}) =  {\bf M}({\bf d^{N_4}})   
{\bf M}({\bf a^{N_1}}){\bf M}({\bf c^{N_3}}){\bf M}({\bf b^{N_2}}).
\]         
Then \(G\) is subrectangular since 
\({\bf M}({\bf a^{N_1}})\) is subrectangular, and \(G\) has 
only finitely many nonzero columns since 
\({\bf M}({\bf b^{N_2}})\) only has finitely many nonzero columns.
Therefore \(G\) can be written
\[
G= 
\left( \begin{array}{cc}
A\;B\; 0\\
0\;\;0\;\; 0\\
C\;D\; 0\\
\end{array}
\right)
\]
where \(A\) is a  \(S_0\times S_0\) matrix with strictly positive 
elements and where  
\(S_0\) is a nonempty finite subset of \(S\), 
\(B\) is an \(S_0\times S_1\)
matrix where  \(S_1 \subset S\) is  finite or empty,
and 
each \(0\) denotes a zero-matrix of appropriate format.

By induction it is straight forward to prove that 
\[G^{n} = 
\left( \begin{array}{cc}
A\\
0\\
C\\
\end{array}
\right)
A^{n-2}
\left( \begin{array}{cc}
A\;B
\;\;0\;
\end{array}
\right)
\]
if \(n \geq 2 \), 

Since \(A\) is a finite-dimensional 
square matrix with
strictly positive elements it follows by 
Perron's theorem (see e.g \cite{Gan65}, vol II, Theorem 8.1),
that there exist a number \(\lambda>0\) and a rank 1 matrix
\({\overline A}\) with strictly positive elements, 
such that 
\[
\lim_{n \rightarrow \infty} ((A^n)_{i,j}/\lambda^n -
({\overline A})_{i,j}) = 0, \;\; \forall i \in S_0, \; j \in S_0.
\]

Finally defining  
\[ W_0= (1/\lambda^2) 
\left( \begin{array}{cc}
A\\
0\\
C\\
\end{array}
\right)
{\overline A}
\left( \begin{array}{cc}
A\;B\;0\;
\end{array}
\right)
\]
and 
\[
W = W_0/||W_0||
\]
it is an easy matter to show that 
\[
\lim_{n \rightarrow \infty} 
||G^{n}/|| G^{n}||-W|| = 0
\]
from which obviously Condition B1 follows.
\(\;\Box\)

\section{ Random walk examples.}

Let \(S = \{0,1,2,...\}\), let
\(S_{odd} = \{ i \in S, \; i \;\; odd\}\) and
\(S_{even} = \{ i \in S, \; i \;\; even \}.\)
For \(i=1,2,... \), let \(a_i,\;
b_i \;\)
 and \(c_i\)
be positive numbers
satisfying
\[ a_i + b_i + c_i=1\]
and let \(b_0\) and \(c_0\) be positive numbers satisfying 
\[b_0+c_0 = 1.\]

Let 
\(P\in PM(S)\) be defined such that 
\[
(P)_{i,i} =b_i, i \in S,
\]
\[\;(P)_{i,i+1}=c_i, \; i \in S  
\] 
and
\[(P)_{i,i-1}= a_i, \; i =1,2,...\;.
\]

Let \({\cal M} =\{M(1),M(2)\}\) be a partition of \(P\) such 
that if \(i\) is odd then the \(i-th\) column of \(P\) is equal to the
\(i-th\) column of \(M(1)\)
and if \(i\) is even then the \(i-th\) column of \(P\) is equal to the
\(i-th\) column of \(M(2)\).

\begin{theorem}\label{thm10.1}
Let the tr.p.m \(P\) and the partition \({\cal M}\) be defined
as above,
and suppose also that 
\begin{equation}\label{T10eq1}
\sum_{n=1}^{\infty} \prod_{i=1}^n c_{i-1}/a_i <  \infty. 
\end{equation}
A) If \(b_i=b_0, \; \forall i \in S \)
then \({\bf P}_{{\cal M}}\) is asymptotically stable.
\newline
B) If there exists \(i_0 \in S\) such that
\(b_{i_0} > \sup\{ b_i: i \in S, \;i \neq i_0 \}\)
then \({\bf P}_{{\cal M}}\) is asymptotically stable.
\end{theorem}

{\em Proof.} From the definition of \(P\) it is clear
that \(P\) is aperiodic and irreducible. That \(P\) is
positively recurrent follows from (\ref{T10eq1}), ( see e.g 
\cite{Pra65}, Problem 18, chapter 2). 
Therefore by Proposition \ref{prop9.1}
and Theorem \ref{thm1.1} it suffices to verify that Condition B1
is satisfied.

We first consider the case when 
\(b_i=b_0, \; \forall i \in S \).
First let us note that 
\(
||e^iM(1)||=b_i=b_0 \;\;if \;\; i \;\;\;  odd, 
\)
and that also 
\(||e^iM(2)||=b_i=b_0 \;\;if \;\; i \;\;\; even.\)
 Since \(M(1)+M(2)=P\), this
implies that 
\(
||e^iM(1)||=1-b_0 \;\;if \;\; i \;\;\;  even, 
\)
and
\(||e^iM(2)||=1-b_0 \;\;if \;\; i \;\;\; odd.\)
Therefore 
\begin{equation}\label{10.eq60}
||e^iM(1)M(2)||=(1-b_0)^2, \; if\;\; i\;even.
\end{equation}
and
\begin{equation}\label{10.eq61}
||e^iM(1)M(2)||=b(1-b_0), \; if\;\; i\;odd.
\end{equation}

Set \(1-b_0=\alpha\) and define
\(M = M(1)M(2)\). Since \((M(2))_{i,j}=0\) if \(j\) is odd it follows
that \((M)_{i,j}=0 \) if \(j\) is odd. 

Now  define the \(S\times S\) matrices \(A\) and \(B\)
by
\[ (A)_{i,j}= (M)_{i,j}, \; i\; is \; even, \;and \; j\;is \; even \]
\[ (A)_{i,j}= 0, \; otherwise \]
and 
\[ B=M-A.\]
It is easily proved that \(BB=0\) and also that \(AB=0\) from which follows
that
\begin{equation}\label{10.eq11}
M^n = (A+B)^n = A^n + BA^{n-1} = MA^{n-1}.
\end{equation}
Since 
\[
||e^iM||=||e^iM(1)M(2)||=(1-b_0)^2= \alpha^2 \; if \;i \;is\; even 
\]
 it
follows that \(\sum_{j \;even} (A)_{i,j} = \alpha^2 \) if \(i \) is
even 
and therefore \(A/\alpha^2\) con be considered as a tr.pr.m on the
even integers. It is easily seen that the induced Markov chain is
irreducible and aperiodic.
Moreover since
\[\sum_{n=0}^\infty \prod_{i=0}^n (M)_{2i,2(i+1)}/(M)_{2(i+1),2i} =\]
\[\sum_{n=0}^\infty \prod_{i=0}^n c_{2i}c_{2i+1}/(a_{2i+2}a_{2i+1})=
\sum_{n=0}^\infty \prod_{i=0}^{2n+1} c_{i}/a_{i+1} < \infty\]
because of (\ref{T10eq1}), it follows that the Markov chain generated
by this
tr.pr.m is also
positively recurrent. Therefore  
there exists 
a probability vector \(q= (q)_i, i \in S\) such that 
\[(q)_{i} = 0, \;\;if \;\;i \;is \; odd, \]
\[
qA/(\alpha^2)=q\]
and if the vector \(y=((y)_i, i \in S)\) is such that \(||y||=1\),
\((y)_i \geq 0, i \in  S\) and \((y)_i = 0, \; if \; i \; odd\) then
\begin{equation}\label{10.eq12}
\lim_{n \rightarrow \infty} ||yA^n/\alpha^{2n} -q||= 0.
\end{equation}
(See e.g \cite{Lin92}, Section 2.1.)

Now suppose 
\[
\alpha \geq 1-\alpha.
\]
From (\ref{10.eq60}) and (\ref{10.eq61}) follows that
\[
||e^iM^n|| = \alpha^{2n}, \; \; i \; is \; even
\]
and 
\[
||e^iM^n|| = (1-\alpha)\alpha^{2n-1}, \; \; i \; is \; odd.
\]
Therefore, when \(\alpha \geq 1-\alpha\) it follows that 
\(||M^n|| = \alpha^{2n}\).

Now define the \(S\times S\) matrix \(W\) by
\[(W)_{i,j}=(q)_j, \; j \in S, \; i \;even ,\]
\[(W)_{i,j}=(1-\alpha) (q)_j/\alpha, \; j \in S, \; i \;odd ,\]
Since \(q\) is a probability vector \(||q|| =1 \) and therefore
\(||W||=1\) since we have assumed that \(\alpha \geq 1-\alpha\).
That \(W \in {\tt W}\) is also clear. 

We now consider \(e^iM^n/||M^n||\) when \(i\) is even. 
Since 
\(||M^n||=\alpha^{2n}\) 
it follows from (\ref{10.eq11}) that
\[
e^iM^n/||M^n|| = e^iM^n/\alpha^{2n} = 
e^iMA^{n-1}/\alpha^{2n} = f^iA^{n-1}/\alpha^{2(n-1)}\]
where \(f^i\) is defined by \(f^i =e^iM/\alpha^2\). 

From the definition of \(f^i\) it is clear that \((f^i)_j\geq 0, \;j \in S\).
Since \(||e_iM||=\alpha^2\) when
\(i\) is even it follows that \(||f^i||=1\) and since
\((M(2))_{i,j}=0\)
if \(j\) is odd it also follows that \((f^i)_j=0\) if \(j\) is odd.
Therefore by (\ref{10.eq12}) we can conclude that
\begin{equation}\label{10.eq14}
\lim_{n \rightarrow \infty} 
||(e^iM^n/||M^n|| - e^iW|| =
\lim_{n \rightarrow \infty} 
||f^iA^{n-1}/(\alpha)^{n-1} - q || = 0
\end{equation}
if \(i\) is even. 

To prove that (\ref{10.eq14}) holds when \(i\) is odd can be done similarly.
Set \(f^i =e^iM/(1-\alpha)\alpha\). Since
\(||e^iM||=(1-\alpha)\alpha\)
when \(i\) is odd
it follows that \(||f^i||=1, (f^i)_j \geq 0, \;j \in S,\; and \; \;(f^i)_j=0, 
\; j\;is\; odd\). 
Hence, when \(i\) is odd, it follows from (\ref{10.eq11}) that
\[
e^iM^n/||M^n|| = e^iM^n/\alpha^{2n} = 
e^iMA^{n-1}/\alpha^{2n} = (1-\alpha)\alpha f^iA^{n-1}/\alpha^{2n}=\]
\[((1-\alpha)/\alpha) f^iA^{n-1}/\alpha^{2(n-1)}.\]
Therefore by (\ref{10.eq12}) we can conclude that
\[
\lim_{n \rightarrow \infty} 
||e^iM^n/||M^n|| - e^iW|| =
\lim_{n \rightarrow \infty} 
((1-\alpha)/\alpha) || f^iA^{n-1}/\alpha^{2(n-1)} - q || = 0
\]
also if \(i\) is odd. Thereby we have verified Condition B1 when
\(\alpha \geq 1-\alpha\). 

Now assume that \(\alpha < 1 - \alpha\). In this case 
\(||M^n||=(1-\alpha)\alpha^{2n-1}\).
We now define 
the \(S\times S\) matrix \(W\) by
\[(W)_{i,j}=\alpha (q)_j/(1-\alpha), \;\; j \in S, \; i \;even ,\]
\[(W)_{i,j}=(q)_j, \;\; j \in S, \; i \;odd .\]
Again it is clear that \(W\) is a rank 1 matrix of norm 1.
By the same reasoning as above we can again conclude
that
\[
\lim_{n \rightarrow \infty} 
||e^iM^n/||M^n|| - e^iW|| = 0,
\;\; \forall \; i \in S.
\]
Thereby we have verified Condition B1 under hypothesis A).

It remains to  consider the case when
there exists \(i_0 \in S\) such that
\(b_{i_0} > \sup\{ b_i: i \in S, \;i \neq i_0 \}\).
Let us first assume that \(i_0 \in S_{odd}\).
Set \[\alpha=\max \{c_{i_0-1},b_{i_0},a_{i_0+1}\},\]
\[
\beta = \max \{b_{i_0}, 1 -b_{i_0}\},
\]
\[\rho = \sup\{ b_i: i \in S, \;i \neq i_0 \}/b_{i_0}\]
and define the vector \(u_0\) by
\[
(u_0)_{i_0-1}=c_{i_0-1}/\alpha,\;\;(u_0)_{i_0}=b_{i_0}/\alpha,\;\;
(u_0)_{i_0+1}=a_{i_0-1}/\alpha,\]
\[
(u_0)_{i}=0,  \;\; if\;\; i \not =\; i_0-1,\;i_0 \;\; or \;\;i_0+1.\]

Furthermore if we define 
the matrix \(D = \{(D)_{i,j} : i \in S, j \in S\}\) by
\[
(D)_{i,i} = (M(1))_{i,i} 
\]
\[
(D)_{i,j} = 0, \;if \;\;i\not = j.
\]
we note that
\[
M(1)^n = M(1)D^{n-1}
\]
and since \(b_{i_0}> b_i\) if \( i \not = i_0\) it follows by elementary
calculations 
that there exists a constant \(C\) independent of \(n\) and \(i
\in S\)  
such that 
\[ 
|(M(1)^n)_{i_0,i_0}/b_{i_0}^n - 1| < C \rho^n
\]
\[
|(M(1)^n)_{i_0-1,i_0}/b_{i_0}^n - c_{i_0-1}/b_{i_0}| < C \rho^n
\]
\[
|(M(1)^n)_{i_0+1,i_0}/b_{i_0}^n - a_{i_0+1}/b_{i_0}| < C \rho^n
\]
and
\[
|(M(1)^n)_{i,j}/b_{i_0}^n| < C \rho^n , \; otherwise.
\] 
Hence,  if we define
the matrix \(W \in {\tt W}\) by
\[
W={u_{0}}^c e^{i_0}
\]
if follows easily that
that there exists a constant \(C_1\) 
independent of \(n\) and \( i \in S\)
such that
\[
|(M(1)^n)_{i,j}-(W)_{i,j}|<C_1 \rho^n, \;\; \forall \;i,j \in S
\]
from which it easily follows that we also have 
\[
\lim_{n \rightarrow \infty} ||M(1)^n/||M(1)^n||-W|| = 0
\]
which of course implies Condition B1.

If instead  there exists an integer  \(i_0 \in S_{even }\) such that
\(b_{i_0} > \sup\{ b_i: i \in S, \;i \neq i_0 \}\)
then by considering 
\(M(2)^n/||M(2)^n||\)
instead of 
\(M(1)^n/||M(1)^n||\) as \(n\) 
tends to infinity, we 
by arguments similar to those 
given above can again prove that Condition B1 is satisfied. 
We omit the details.
\(\;\Box\)

\section{ Exceptional cases.}
In this section
we first present an example from 1974 due to H. Kesten
(\cite{Kes74}), which
shows that the tr.pr.f \({\bf P}_{{\cal M}}\) induced by a partition
of an aperiodic and irreducible tr.pr.m \(P\)
may even  turn out to be 
{\bf periodic}.
We shall then present a theorem
with hypotheses that
guarantee that \({\bf P}_{{\cal M}}\)
is not asymptotically stable
and end this
section with presenting a whole class of tr.pr.ms for which it
is possible to find partitions such that 
\({\bf P}_{{\cal M}}\) is {\bf not} asymptotically stable.

Kesten's example is an extension
of the example in \cite{Kai75} and reads as follows. 
\begin{example}\label{ex11.1}
(\cite{Kes74}) Let \(S=\{1,2,...,8\}\) and define 
\(P \in PM_{ae}(S)\) by
\[
P =
\left( \begin{array}{cc}
\star \;\;0\;\;\,0\;\;\,0\;\,\star\;\;0\;\;\,0\;\;0 \\
\,0\;\, \star\;0\;\;0\;\;\,0\;\;\star\;\;0\;\; 0\\
0 \;\;0\;\;0\;\;\star\;\;0\;\;0\;\;0\;\;\star\\
\;0\;\;\, 0\;\star\;\;0\;\;0\;\;0\;\;\star\;\,0\\
\star \;\,0\;\;\,0\;\;\,0\;\;\,0\;\;0\;\;\,0\;\;\star \\
\,0\;\; \star\;\,0\;\;0\;\;0\;\;0\;\;\star\;\; 0\\
0 \;\;0\;\;\,0\;\;\star\;\;\star\;\;0\;\;0\;\;0\\
\,0\;\; 0\;\;\star\;\;0\;\;0\;\;\star\;\;0\;\;0\\
\end{array}
\right),
\]
where each \(\star\) denotes the value 1/2.
Let \(A=\{a,b\}\), let \(g:S\rightarrow A\) be a ``lumping''
function defined by \(g(i)=a, \; i=1,2,3,4\), 
\(g(i)=b, \; i=5,6,7,8\), and let \({\cal M}\) be the
partition determined by the ``lumping'' function \(g\). 
(See Example \ref{ex1.1}.)
\end{example} 

That \(P\) is aperiodic and irreducible is easily seen. It is also not
difficult to verify that if for example \(x_0\in K\) is such that
\[
(x_0)_i=0, i=5,6,7,8  \;\;and \;\;
0< (x_0)_1=(x_0)_3 < (x_0)_2= (x_0)_4 
\]
then the Markov chain generated by \({\bf P}_{{\cal M}}(x_0, \cdot) \) 
is a {\bf periodic} Markov chain taking
its values in a subset of \(K\)
consisting of just 8 elements with coordinates  depending on \(x_0\).

We shall next state a theorem with hypotheses that guarantee that 
\({\bf P}_{{\cal M}} \) is {\em not asymptotically stable}.

To state the theorem we need
two further  notations, \(K(x,{\cal M})\) and
\(K_{S'}\). Let \(S\) be a denumerable set and let
\({\cal M}= \{M(w): w \in {\cal W}\}\in G(S)\). 
Recall that \({\cal W}_{{\cal M}}(x)\)
is defined by (\ref{eq1.5.1}).
For each \(x \in K\) we define
\[
K(x,{\cal M}) = \{ y \in K: y= xM(w)/||xM(w)|| 
\;\;some\;\;  w \in {\cal W}_{{\cal M}}(x)\}.\]
Next let \(S' \subset S\). We define
\[
K_{S'}= \{x \in K : (x)_i= 0, \;i \not \in S'\}.\]
\begin{theorem}\label{thm11.1}
Let \(S\) be a denumerable set,  
let \({\cal M}=\{M(w) \,: \, w \in {\cal W}\} 
\in G(S)\) 
and let \({\bf P}_{{\cal M}}\) be the tr.pr.f induced by  
\({\cal M}\).
Suppose that there exists a subset \(S' \subset S\) consisting of
at least two elements, such that
\newline
1) for every \(x \in K_{S'}\) the set 
\(
\cup_{n=1}^\infty K(x,{\cal M}^n)
\)
consists of isolated points,
\newline
2) if both \(x\) and \(y \) are in \(K_{S'}\) 
then
\(
{\cal W}_{{\cal M}^n}(x)= {\cal W}_{{\cal M}^n}(y), \;\;n=1,2,....
\)
\newline
and 
\newline 
3) if \(x\) and \(y \) in \(K_{S'}\), \(n \geq 1\) and  
\({\bf w^{n}} \in {\cal W}_{{\cal M}^n}(x)\) then 
\[
||(x{\bf M}({\bf w^{n}})/||x{\bf M}({\bf w^{n}})|| -
y{\bf M}({\bf w^{n}})/||y{\bf M}({\bf w^{n}})||)|| =||x-y||.
\]

If such a subset \(S'\) exists then
\({\bf P}_{{\cal M}}\) is {\bf not} asymptotically stable.
\end{theorem}
\begin{rem} If we could prove that the hypotheses of 
Theorem \ref{thm11.1} are also
necessary in order for 
\({\bf P}_{{\cal M}}\) to be a tr.pr.f which
is  {\em not asymptotically stable} when \(P \in PM_{ae}(S)\), we would
have 
a rather  easily checked criterion for deciding whether
a tr.pr.f 
\({\bf P}_{{\cal M}}\)
induced by a partition \({\cal M}\) 
is asymptotically stable or
not.
\end{rem}
\begin{rem} It is easy to check that Kesten's example satisfies 
the hypotheses of Theorem \ref{thm11.1}.
\end{rem}
{\em Proof of Theorem \ref{thm11.1}.} 
Let \(S' \subset S\) be as in the
hypotheses of the theorem. Let \(x \in  K_{S'}\)
and set 
\[K'(x)=\cup_{n=1}^\infty K(x,{\cal M}^n).\]
Because of hypothesis 1), the set \(K'(x)\) consists of isolated
points,
and because of hypothesis 3) it is not difficult to convince oneself
that \(K'(x)\) must contain at least two points.
Therefore
\[
\epsilon_0= \inf \{||z_1 - z_2||: z_1, z_2 \in K'(x),\;z_1 \neq z_2 \} 
> 0. 
\]
Since \(x \in K_{S'}\) and \(S'\) consists of at least 
two elements,  we can find an element \(y \in 
K_{S'}\) such that \(||x-y||=\epsilon_0/2.\)

Now let \(n\) denote an arbitrary positive integer, let
\[
K_1= \{xM({\bf w^n})/||xM({\bf w^n})||: 
{\bf w^{n}}\in {\cal W}_{{\cal M}^n}(x)\}
\]
and
\[
K_2= \{yM({\bf w^n})/||yM({\bf w^n})||: 
{\bf w^{n}}\in {\cal W}_{{\cal M}^n}(x)\}.
\]

Since the points in \(K'(x)\) are isolated points, it is clear that
\(K_1\) is the support of \({\bf P}^n_{{\cal M}}(x,\cdot)\)  
and from hypotheses 2) 
that also \(K_2\) is the support of 
\({\bf P}^n_{{\cal M}}(y,\cdot)\).
Since 
\[\inf \{\delta(z,K_1): z \in K_2\} = \epsilon_0/2 \]
because of hypothesis 3), it therefore follows that 
the Kantorovich distance 
\[
d_K({\bf P}^n_{{\cal M}}(x,\cdot),{\bf P}^n_{{\cal M}}(y,\cdot)) \geq 
\epsilon_0/2 > 0
\]
and since \(n\) was an arbitrary integer 
\[
d_K({\bf P}^n_{{\cal M}}(x,\cdot),{\bf P}^n_{{\cal M}}(y,\cdot)) \geq 
\epsilon_0/2 
\]
for \(n=1,2,...\) which implies that the tr.pr.f
\(\;{\bf P}^n_{{\cal M}}(\cdot,\cdot)\) can not be asymptotically stable.
\(\;\Box\)

We shall next describe a family of 
tr.pr.ms 
for which  one, for each matrix belonging to the family, can find
a partition such that  the induced
tr.pr.f
is {\em not} asymptotically stable.

Let \({\cal I}\) denote
a denumerable set, let \(d\geq 2\) be an integer,  
set \(I_d=\{1,2,...,d\}\), and define
the 
set \(S\) 
as the Cartesian product of 
\({\cal I}\) and \(I_d\), that is 
\[
S = \{(i,j), i \in {\cal I}, j \in I_d \}.
\]
Let  
\(
A \in PM({\cal I}),
\)
and let 
\({\cal M}= \{ M(w):\; w \in  {\cal W}\}\)
be a partition of \(A\).

Next  let 
\(Perm(d)\) 
denote the set of \(d\times d\) 
{\em permutation matrices}.
For each \(w \in {\cal W}\) and each \((i,k)\in {\cal I}\times
{\cal I}\) 
we now associate a matrix \(Q(i,k,w) \in 
Perm(d)
\). We write
\[
{\cal Q}_{{\cal I},{\cal W}} = \{ Q(i,k,w):\; (i,j) \in {\cal I}\times
{\cal I},\;w \in {\cal W}\}.
\] 
We 
define the set
\[
{\cal M}'=\{M'(w): w \in {\cal W}\}
\]
of \(S\times S\) matrices 
by
\begin{equation}\label{U8eq12}
(M'(w))_{(i,j),(k,m)}=(M(w))_{i,k}\cdot (Q(i,k,w))_{j,m},
\end{equation}
and define the \(S\times S\) matrix \(P\) by
\[
(P)_{(i,j),(k,m)}=\sum_{w \in {\cal W}}
(M'(w))_{(i,j),(k,m)}.
\]
It is easily verified that \(P \in PM(S)\) and that 
\(
{\cal M}'\in G(S,P) \).
We call \(P\) the tr.pr.m {\em generated} by \(A\) and 
\({\cal Q}_{{\cal I},{\cal W}}\).

Next suppose that the partition 
\({\cal M}=\{M(w):w \in{\cal W}\}\) of \(A\) 
is such that 
\begin{equation}\label{S8eq100}
(M(w))_{i,k} > 0 \;\; \Rightarrow  \;\;(M(w))_{i,k_1} = 0, \;
if \;\;k_1 \neq k, \;\;\forall M(w)\in {\cal M}.   
\end{equation}

\begin{proposition}\label{prop11.2}
Let \({\cal I}\) denote
a denumerable set, 
let \(A\in PM({\cal I})\) and let 
\({\cal M}=\{M(w) : w \in {\cal W}\}\)
be a partition of \(A\) that satisfies (\ref{S8eq100}). 
Let \(d\geq 2\), set \(I_d=\{1,2,...,d\}\), and let
\[
{\cal Q}_{{\cal I},{\cal W}} = \{ Q(i,k,w):\; (i,j) \in {\cal I}\times
{\cal I},\;w \in {\cal W}\}.
\] 
Let 
\(S = \{(i,j), i \in {\cal I}, j \in I_d \}\), 
let
\(P \in PM(S)\) be the tr.pr.m generated by \(A\)
and 
\({\cal Q}_{{\cal I},{\cal W}}\)
and let \({\cal M}'\) be the partition of \(P\) defined 
by (\ref{U8eq12}).

Then 
\({\cal M}'\) satisfies
the hypotheses of Theorem 11.1.
\end{proposition}
{\em Proof.}
The proof is based on the following observation. 
For \(i \in {\cal I} \), let \(S'_i =\{(i,j),j=1,2,...d\}\). 
Let \(x \in K_{S'_i}\) 
and suppose 
\(||xM'(w)|| > 0\). From (\ref{S8eq100}) and (\ref{U8eq12}) follows that 
\(xM'(w)/||xM'(w)|| \in K_{S'_k}\) where thus \(k\) is such that
\((M(w))_{i,k} > 0\). Furthermore, if we let \(x'\) denote the 
{\em d-dimensional} vector defined by 
\((x')_j= (x)_{i,j}, \; j=1,2,...,d\),
set \(z=xM'(w)/||xM'(w)||\) 
and let \(z'\) denote the {\em d-dimensional} vector 
defined by
\((z')_j =(z)_{k,j},\;\; j=1,2,...,d,\;\) then
\(z'= x'Q(i,k,w)\).
Since \(i\) was arbitrary it now easily follows that 
 the hypotheses of Theorem 11.1 are fulfilled. \(\;\bullet \)

It is easy to show that both the example in \cite{Kai75} and 
Kesten's example can be put
into
the framework of the class just described. (Choose 
\({\cal I}=\{1,2\}, \;d=2\) for the example in \cite{Kai75}
and
\({\cal I}=\{1,2\}, \;d=4\) for Kesten's example.)
When doing this for
Kesten's example it turns out that all permutation  matrices 
will be {\em odd} and that explains why the induced tr.pr.f 
 gives rise to a {\bf periodic} 
Markov chain for most initial distributions.
\begin{conjecture}\label{conj11.3}
 If \(S\) is denumerable, \(P \in PM_{ae}(S)\),
\({\cal M} \in G(S,P)\) and \({\bf P}_{\cal M}\) is 
{\bf not} asymptotically
stable, then \(P\) and \({\cal M}\) can be represented as in 
Proposition 11.1.
\end{conjecture}
\begin{conjecture}\label{conj11.4}
Suppose  \(S\) is a finite set with size equal to
a {\bf prime number}, let  \(P \in PM_{ae}(S)\), 
and let \({\cal M}\) be a partition of \(P\) determined by
a ``lumping'' function. Then \({\bf P}_{\cal M}\) has
a unique invariant measure.
\end{conjecture}
In other words: {\em We believe that Blackwell's conjecture 
(see \cite{Bla57}, page 19) is true if we
add the hypothesis ``the number of states is a prime number''.}

Before finishing this section let us mention that if 
\(S\) is a finite set of size \(d\geq 2\), and the tr.pr.m 
\(P \in PM(S)\) is {\bf doubly stochastic}, then it is always 
possible to find a partition of \(P\) such that
the hypotheses of Theorem 11.1  are fulfilled since
every doubly stochastic tr.pr.m can be written as 
a weighted sum of permutation matrices. Obviously 
a doubly stochastic tr.pr.m \(P\) belongs to the class considered in 
Proposition 11.1 since we can choose the set \({\cal I} =\{1\},\)
the matrix \(A\) equal to 1, and the partition of 
\(A\) equal to the
weights used in the representation of \(P\) as a weighted
sum of permutation matrices.

\section{Convex functions and barycenters.}

As usual, let \(S\) be a denumerable set and let \(K\)
denote the set of probability vectors on \(S\). Clearly
\(K\) is a convex set. Let  
\(C_{convex}[K]\) denote the subset of \(C[K]\) consisting of
convex functions.
In this section we shall present some simple inequalities 
for functions in \(C_{convex}[K]\).

 Let 
\(C'_{convex}[K]\) 
denote the subset of \(C[K]\) that 
can be obtained as 
\[ u = \sup \{v_n : n \in {\cal N}\}
\]
where \({\cal N}\) is an arbitrary  index set, and each \(v_n\) is
an affine function on \(K\) such that 
\[
v_n(x) = xa_n^c+b_n  
\]
where \(a_n\in l^{\infty}(S)\)
and each \(b_n\) is a real number.
By using the fact that it is possible to
find a separating affine plane between the epigraph of a convex
function on
\(K\) and a point outside the epigraph (because of Hahn-Banach's
theorem)
 it is not difficult to prove that
 \(C_{convex}[K]=C'_{convex}[K]\).

Next, for each \(q \in K\) we define the subset
\({\cal P}_d(K|q)\) of \({\cal P}(K|q)\) as 
the set consisting of all
measures \(\nu\)
such that  
\( \nu = \sum_{k=1}^{\infty} \alpha_k \delta_{x_k} \)
where
\(x_k,\; k=1,2,...\;\in K\) 
and
\(\sum_{k=1}^\infty \alpha_k = 1 \).
We define  \(\psi_q \in {\cal P}(K|q)\) by 
\(\psi_q(e^i)=q_i, \;i \in S\).

The following proposition is a simple consequence of convexity.
\begin{proposition}\label{prop12.1} 
Let \(\mu \in {\cal P}_d(K|q)\) and
\(u \in C_{convex}[K]\). 
Then
\begin{equation}
\langle u,\delta_q \rangle  
 \leq \langle u,\mu\rangle \leq \langle u,\psi_q\rangle. 
\end{equation}
\end{proposition}
{\em Proof.} Let \(u \in C_{convex}[K]\) and
\( \mu = \sum_{k=1}^{\infty} \alpha_k \delta_{x_k}\in {\cal
  P}_d(K|q)\). Then
\[
\langle u,\delta_q \rangle =u(q)\leq
\sum_{k=1}^{\infty} \alpha_k u(x_k)= 
\langle u,\mu\rangle \leq 
\sum_{k=1}^{\infty}\alpha_k (\sum_{i=1}^{\infty}(x_k)_i u(e^i))=\]
\[\sum_{i=1}^{\infty}(\sum_{k=1}^{\infty}\alpha_k (x_k)_i u(e^i))=
\sum_{i=1}^{\infty}q_i u(e^i) = \langle u, \psi_q \rangle. \;\;\Box\]

Also the next proposition is easy to prove if one uses the fact that 
if \(u \in C_{convex}[K]\) then \(u \in C'_{convex}[K]\)
and the fact that 
for any index set \({\cal N}\) 
\[\sup\{f_n(x)+g_n(x): n\in {\cal N}\} 
\leq \sup \{f_n(x): n\in {\cal N}\} +
\sup\{ g_n(x): n\in {\cal N}\}.\]
We omit the details.  The result is essentially due  to 
H. Kunita. (See \cite{Kun71}, Lemma 3.2.)
\begin{proposition}\label{prop12.2}
Let \(S\) be a denumerable set  and let \({\cal M} \in {\cal G}(S).\)
Then 
\(u \in 
C_{convex}[K]  \Rightarrow T_{{\cal M}}u \in 
C_{convex}[K].
\)
\end{proposition}

Our next proposition 
follows easily by using Theorem \ref{thm6.2}, formula (\ref{I1eq1.5.3}),
Proposition \ref{prop12.1} and 
Proposition \ref{prop12.2}

\begin{proposition}\label{prop12.3}
Let \(S\) be a denumerable set, let \(P_1, P_2 \in PM(S)\) let
\({\cal M}_1 \in G(S,P_1)\), let 
\({\cal M}_2 \in {\cal G}(S,P_2)\)
and suppose that 
\(\pi \in K\) is such that 
\[\pi=\pi P_1=\pi P_2.\]  
Then, if \(u \in C_{convex}[K] \)
\[
\langle u,\delta_\pi \rangle \leq 
\langle u,{\breve P}_{{\cal M}_2}\delta_\pi \rangle \leq 
\langle u,{\breve P}_{{\cal M}_2}
{\breve P}_{{\cal M}_1}\delta_\pi\rangle \leq
\langle u,{\breve P}_{{\cal M}_2}
{\breve P}_{{\cal M}_1}\psi_\pi\rangle \leq
\langle u,{\breve P}_{{\cal M}_2}\psi_\pi\rangle \leq 
\langle u,\psi_\pi\rangle.
\]
\end{proposition}
{\em Proof.} 
The first inequality follows from Theorem \ref{thm6.2},  
and Proposition \ref{prop12.1}.
The second inequality follows by using 
(\ref{I1eq1.5.3}) and the first inequality of 
Proposition \ref{prop12.3} which we just proved.
The third inequality follows by using (\ref{I1eq1.5.3}), Proposition
\ref{prop12.2} 
and Proposition \ref{prop12.1}.
The last inequality 
follows from Proposition \ref{prop12.1} and Theorem \ref{thm6.2}, 
and finally the forth inequality follows
by using 
(\ref{I1eq1.5.3}) together  the last inequality of 
Proposition \ref{prop12.3}  which we just proved.
\(\; \Box\)

As an immediate corollary we obtain 
\begin{corollary}\label{cor12.4}
Let \(S\)
be a denumerable space, 
let
\(P \in PM(S) \), let
\(\pi \in K\)
satisfy
\(\;\pi =\pi P,\;\)  
let
\({\cal M}\in G(S,P)\)
and let \(u \in C_{convex}[K]\). Then, for \(n=1,2,...,\;\;\)
\[
\langle u, {\breve P}_{{\cal M}}^n\delta_\pi \rangle
\leq
\langle u, {\breve P}_{{\cal M}}^{n+1}\delta_\pi \rangle
\leq
\langle u, {\breve P}_{{\cal M}}^{n+1}\psi_\pi \rangle
\leq
\langle u, {\breve P}_{{\cal M}}^{n}\psi_\pi \rangle.
\]
\end{corollary}
{\em Proof.} Follows from Proposition \ref{prop12.3} 
and Corollary \ref{cor3.3} \(\;\Box.\)

\section{ Blackwell's entropy formula.}
We now return to the same set-up as the one we started with in 
subsection 1.1. Thus,
let \(S\) denote a denumerable
set, \(\{X_n, n=0,1,2,...\}\) be an  aperiodic, 
positively recurrent  Markov chain with  
\(S\) as state space,
let \({\cal W}\) denote an ``observation space'',
let  \(g:S \rightarrow {\cal W}\) denote a "lumping" function of 
the state space \(S\) and define 
\(Y_n = g(X_n)\).

In \cite{Bla57} D. Blackwell
presented an integral formula
for the entropy rate of the \(Y_n-process\) when the \(X_n-process\) is
a stationary process and \(S\) is finite.
 The purpose of this final section is to
generalize
Blackwell's result.

Let  
\({\cal M}=\{M(w): w \in {\cal W}\}\in G(S,P)\) be determined
by the "lumping" function \(g\), and assume that
\begin{equation}\label{I2eq21}
\sup_{x \in K} \sum_{w \in {\cal W}_{\cal M}(x)}
- (\ln||xM(w)||)\cdot ||xM(w)|| <  \infty.
\end{equation}
Let \(\pi \in K\) satisfy \(\pi P=\pi\) and let 
 \(h:[0, 1] \rightarrow [0, 1/(e\cdot \ln(2))] \) be defined by
\[
h(t)=-t\ln(t)/\ln(2), \;\;if\;\;0 < t \leq 1  \;\;\;\;\;\;\;\;\;\;
and \;\;\;h(0)=0.
\]
 For \(x\in K\) we define
\[
H^n(Y;x) = 
\sum_{{\bf w^n} \in {\cal W}^n}h(||x {\bf M}({\bf w^n})||)
\]
and 
\[
H_R^n(Y;x) = H^{n+1}(Y;x) - H^n(Y;x).  
\]

\begin{theorem} \label{thm13.1} 
Let \(S\) denote a denumerable set, let
\(P\in PM_{ae}(S)\), let \(g\) be a ``lumping'' function on \(S\),
let \({\cal M}=\{M(w): w \in {\cal W}\}\in G(S,P)\) be determined
by the "lumping" function \(g\) and let
\({\bf P}_{{\cal M}}\) be the tr.pr.f induced by \({\cal M}\).
Suppose also that Condition B is satisfied and that
(\ref{I2eq21}) holds.
Then,  
\newline
a) 
for every \(x \in K\) and 
\(n=1,2,...\)
\[
H_R^n(Y;x) = \sum_{w \in {\cal W}} \int_{K} h(||yM(w)||)
{\bf P}_{{\cal M}}^n(x,dy)
\]
b) for \(n=1,2,...\)
\[
\int_{K} H_R^n(Y;x)\psi_{\pi}(dx) \leq
\int_{K} H_R^{n+1}(Y;x)\psi_{\pi}(dx) \leq
\]
\[
H_R^{n+1}(Y;\pi) \leq
H_R^{n}(Y;\pi)
\]
c) for every \(x\in K\)
 \[
\lim_{n \rightarrow \infty}
H_R^n(Y;x) = 
\sum_{w \in {\cal W}} \int_{K}h(||yM(w)||)\mu(dy)
\]
where \(\mu\) is the unique invariant measure of 
the tr.pr.f \({\bf P}_{{\cal M}}\).
\end{theorem}
{\em Proof.} From the definition of 
 \(H^n(Y;x)\), (\ref{I2eq21}) and (\ref{G2eq11}) we find that 
\[
H^{n+1}(Y;x) =\]
\[(-1/\ln(2)) \sum_{w_{m+1}\in {\cal W}}
\sum_{{{\bf w^n}\in {\cal W}_{{\cal M}^n}^n}(x)} 
(\ln(||(x{\bf M}({\bf w^n})/||x{\bf M}({\bf w^n})||)M(w_{m+1})||)\cdot
\]
\[(||(x{\bf M}({\bf w^n})/||x{\bf M}({\bf w^n})||)M(w_{n+1})||)\cdot
||x{\bf M}({\bf w^n})|| +\]
\[(-1/\ln(2))\sum_{w_{m+1}\in {\cal W}}
\sum_{{{\bf w^n}\in {\cal W}_{{\cal M}^n}^n}(x)} 
\ln(||x{\bf M}({\bf w^n})||)
||x{\bf M}({\bf w^n})M(w_{n+1})||=
\]
\[
\sum_{w \in {\cal W}} \int_{K}h(||yM(w)||){\bf P}_{{\cal M}}^n(x,dy) +
 H^n(Y;x) 
\]
which proves a).
Next, using the fact that the function  
\(f_w: K \rightarrow {\tt R}\) defined by \(f_w(y)= - h(||yM(w)||)\) 
belongs to \(C_{convex}[K]\) for every \(w \in {\cal W}\) 
the assertions in b) and c) follow from Corollary \ref{cor12.4} and
Theorem \ref{thm1.1} respectively together with (\ref{I2eq21}) an inequality
which guarantees
that the sums involved are finite. \(\;\Box\)

\acknowledgements {\rm
I am deeply grateful to my brother Sten Kaijser for
always giving me the opportunity 
to discuss my work with him. 
I also want to thank Arne Enqvist, Torkel Erhardsson   and 
 Robert Forchheimer for valuable and constructive comments.
This work has  been supported by the Institute
of Applied Mathematics, 
Academy of Mathematics and Systems Science, Chinese Academy of
Science.}

\end{document}